% bbk18.tex_Resubmit version 
% TK, 11 April 2005

\magnification 1200
\hsize=17 truecm  
\vsize=8.5  truein
\baselineskip=14 truept

\def\bigb{\bigbreak\noindent}         
\def\med{\medbreak\noindent}
\def\small{\smallskip\noindent}
\def\nl{\hfil\break}

\def\today{\noindent\number\day
\space\ifcase\month\or
  January\or February\or March\or April\or May\or June\or
  July\or August\or September\or October\or November\or December\fi
  \space\number\year}

\def\cl {\centerline}
\def\i {\item}

\input amssym.def
\input amssym.tex
\input epsf

 \def\bE {{\Bbb E}} 
\def\bR {{\Bbb R}} \def\bN {{\Bbb N}} 
\def\R {{\Bbb R}}  
  
  \def\sC {{\cal C}}
 \def\sE {{\cal E}} \def\sF {{\cal F}}
  
  \def\sL {{\cal L}}

\def\q{\quad} \def\qq{\qquad}

\def \half {{\textstyle {1 \over 2}}}

\def\fract#1#2{{\textstyle {#1\over #2}}}

\def\qed {\hfill$\square$\par}

\def\PHI{{\mathop {{\rm PHI}}}}
\def\CS#1{\mathop {{\rm CS}(#1)}}
\def\betaloc{{\bar \beta}}

\def\til-1{{\,\buildrel -1 \over \sim\,}}

\def\lam{{\lambda}}
\def\th{{\theta}}
\def\al{{\alpha}}
\def\Gam{\Gamma}

\def\grad{{\nabla}}
\def\proof{{\medskip\noindent {\it Proof. }}}
\def\longproof#1{{\medskip\noindent {\bf Proof #1.}}}
\def\qed{{\hfill $\square$ \bigskip}}
\def\subsec#1{{\bigbreak\noindent \bf{#1}.}}

\def\section#1#2{{\bigbreak\bigskip \centerline{\bf #1. #2}\bigskip}}

\def\cite#1{{[#1]}}

\def\eps{\varepsilon}
\def\vp{\varphi}

\def\pd {\partial}

\def\q {\quad} \def\qq {\qquad}
\def\frac#1#2{{#1\over #2}}

\def\del{{\partial}}
\def\wt{\widetilde}
\def\ol{\overline}

\def\wh{\widehat}

\def\ni{\noindent}

\def\ms{\medskip}
\def\bs{\bigskip}
\def\cl#1{\centerline{#1}}

\def\ftN{\rl\kern-0.13em\rN}

\def\E{{\Bbb E}}

\def\dm {d\mu} \def\dmp{} 

\def\tlint{{- \kern-0.85em \int \kern-0.2em}}  % for textstyle
\def\dlint{{- \kern-1.05em \int \kern-0.4em}}  % for displays

%%%%Def from Martin's survey
\def\Lam {\Lambda}
\def\dg{d_G}
\def\wvp{\wh\vp}
\def\gam{\gamma}

\def\loc{{\rm loc}}
\def\Reff{R_{\rm eff}}
\def\PHI{\mathop{\hbox{PHI}}}
\def\EHI{\mathop{\hbox{EHI}}}
\def\VD{\mathop{\hbox{VD}}}
\def\PI{\mathop{\hbox{PI}}}
\def\CS{\mathop{\hbox{CS}}}
\def\CC{\mathop{\hbox{CC}}}
\def\HK{\mathop{\hbox{HK}}}
\def\RES{\mathop{\hbox{RES}}}

%%%%New Def 
\def\es{\hbox{ess sup}}
\def\ei{\hbox{ess inf}}
\def\EE{\mathop{\hbox{E}}}

%%% Start of MSS

\font\titlefont=cmbx10 scaled\magstep2
                
\centerline {\titlefont Stability of parabolic
Harnack inequalities} 
\centerline {\titlefont on metric measure spaces }
\medskip

\vskip .5truein

\centerline {Martin T. Barlow\footnote{$^1$}{Research partially
 supported by a grant from NSERC (Canada).},
\qquad  Richard F. Bass\footnote{$^2$}{Research partially supported by
NSF Grant DMS-0244737.} \qquad and \qquad 
Takashi Kumagai\footnote{$^3$}{Research partially supported by the
Program for Overseas Researchers of the Japan Ministry of Education and by the 
Grant-in-Aid for Scientific Research for Young Scientists (B) 16740052.}}

\vskip 0.6 truein

\ni{\bf Abstract.} Let $(X,d,\mu)$ be a metric measure space with a 
local regular Dirichlet form. We give necessary and sufficient
conditions for a parabolic Harnack inequality with global space-time
scaling exponent $\beta\ge 2$ to hold.  We show that this parabolic
Harnack inequality is stable under rough isometries. As a consequence,
once such a Harnack inequality is established on a metric measure
space, then it holds for any uniformly elliptic operator in divergence
form on a manifold naturally defined from the graph approximation of
the space.

\bs
\bs
\ni{\bf Key words and Phrases.} Harnack inequality, volume doubling, Green functions, 
Poincar\'e inequality, Sobolev inequality, rough isometry,  
anomalous diffusion

\bs
\bs
\ni{\bf Short title.} Harnack inequalities on metric spaces                     

\bigskip
\bigskip

\ni{\bf 2000 Mathematics Subject Classifications.} 
Primary: 60J35; Secondary: 31B05, 31C25

\vfill\eject
  
\subsec{1. Introduction}
Let $(X,d,\mu)$ be a metric measure space with a local regular
Dirichlet form $(\sE,\sF)$.  (See Section 2 for definitions of the
terms used in the introduction.)  Assume that the metric is geodesic
and for simplicity assume that $X$ has infinite diameter. For a
typical example let $X$ be a complete non-compact Riemannian manifold
with Riemannian metric, Riemannian measure and $\sE(f,f)=\int_X|\nabla
f|^2d\mu$. In this paper, we give necessary and sufficient conditions
for a parabolic Harnack inequality with `anomalous'  space-time scaling
to hold, and we show the stability of this
parabolic Harnack inequality under certain transformations.

For $\betaloc, \beta\ge 2$, 
let $\Psi (s)=s^{\betaloc}1_{\{s\le 1\}}+s^{\beta}1_{\{s> 1\}}$.
For $R>0$ let $T=\Psi(R)$, 
$Q_- = (T,2T)\times B(x_0,R)$ and $Q_+ = (3T,4T)\times B(x_0,R)$. 
Let $\Delta$ be the self-adjoint operator corresponding to $(\sE,\sF)$. 
$X$ satisfies the {\sl parabolic Harnack inequality} $\PHI(\Psi)$,
if there  exists a constant $c_{1}$ such that for any $x_0\in X$ and $R>0$, 
if $u=u(t,x)$ is a non-negative solution of the heat equation
$\pd_t u = \Delta u$ in $(0, 4T)\times B(x_0, 2R)$, then 
$$ \sup_{Q_-} u \le c_1 \inf_{Q_+} u.$$
If $\betaloc=\beta$ we sometimes write $\PHI(\beta)$ for  $\PHI(\Psi)$.
It is known that such a Harnack inequality is strongly related to detailed
estimates of the heat kernel, i.e.,  the fundamental solution of the heat equation. 

When $\betaloc=\beta= 2$, the study of the parabolic 
Harnack inequality has a long history
(see [Da, SC2] etc. for details).
For any divergence operator ${\cal L}=\sum_{i,j=1}^n\frac{\partial}{\partial x_i}
(a_{ij}(x)\frac{\partial}{\partial x_j})$ on $\bR^n$ which is 
uniformly elliptic, Aronson ([A]) proved that the heat kernel $p_t(x,y)$ satisfies 
$$\frac{c_1}{\mu (B(x,t^{1/2}))}\exp \Big(-\frac{d(x,y)^2}{c_1t}\Big)\le p_t(x,y)\le 
\frac{c_2}{\mu (B(x,t^{1/2}))}\exp \Big(-\frac{d(x,y)^2}{c_2t}\Big),\eqno (1.1)
$$
where $\mu$ is Lebesgue measure 
(so $\mu (B(x,t^{1/2}))^{-1}= c t^{-n/2}$).
It is not hard to derive the parabolic Harnack inequality from (1.1).
Similar results hold 
in the field of global analysis on manifolds. Let $\Delta$ 
be the Laplace-Beltrami operator on a complete Riemannian manifold $M$ with 
Riemannian metric $d$ and with Riemannian measure $\mu$. Li-Yau ([LiY]) proved 
that if $M$ has non-negative Ricci curvature, then the heat kernel $p_t(x,y)$ 
satisfies (1.1). 
A few years later, Grigor'yan ([Gr1]) and Saloff-Coste ([SC1])  
refined this result and proved, in conjunction with results 
by Fabes-Stroock ([FS]) and Kusuoka-Stroock ([KS]), that (1.1) 
is equivalent to a volume doubling condition $(\VD)$ 
plus a Poincar\'e inequality $(\PI(2))$. 
These techniques were then extended to Dirichlet forms on metric spaces 
in [BM, St1, St2] and to 
graphs in [Del]. The origin of the ideas and techniques used in this 
field go back to Nash ([N]) and Moser ([M1, M2, M3]).

Examples where $\betaloc=\beta>2$ are given by fractals.
The mathematical study of stochastic processes and the corresponding operators on fractals 
(see, for instance, [B1, Ki]) has shown that on many  
`regular' fractals, there are naturally defined Dirichlet forms whose heat kernels
with respect to the Hausdorff measure $\mu$ satisfy 
$$ \eqalignno{
\frac{c_1}{\mu (B(x,t^{1/\beta}))}
\exp \Big(-\Big(\frac{d(x,y)^\beta}{c_1t}\Big)^{\frac 1{\beta-1}}& \Big)
\le p_t(x,y) \cr
&\le \frac{c_2}{\mu (B(x,t^{1/\beta}))}
\exp \Big(-\Big(\frac{d(x,y)^\beta}{c_2t}\Big)^{\frac1{\beta-1}}\Big),
 &(1.2)}$$
(see [BP, BB1, BB2, FHK, Kum1] etc.).
The techniques of [Gr1] and [SC1] do not apply to
these spaces. The main obstacle is that the Moser iteration argument
used in these papers needs the
existence of sufficiently many cut-off functions $\vp$
with approximately minimal energy 
such that the $L^2$ and $L^\infty$ norms of $\grad \vp$, suitably
normalized, are comparable. Functions of this kind exist only if $\betaloc=\beta=2$.

Recently, in [BB5], two of the authors established equivalent stable conditions
for (1.2) (equivalently,  to $\PHI(\beta)$) on graphs. The key idea was to 
introduce a new inequality $\CS(\beta)$, called a cut-off Sobolev inequality, 
which implies the existence of enough `low energy' cut-off functions on the space.
This was quite recently extended to the manifold setting by one of the authors in [B2].
Note that when the process is strongly recurrent, there are simpler equivalent stable 
conditions (see [BCK, Kum3]), given in terms of electrical resistance. 
 
\ms ~~~~~\epsfysize 1.8 in   \epsffile{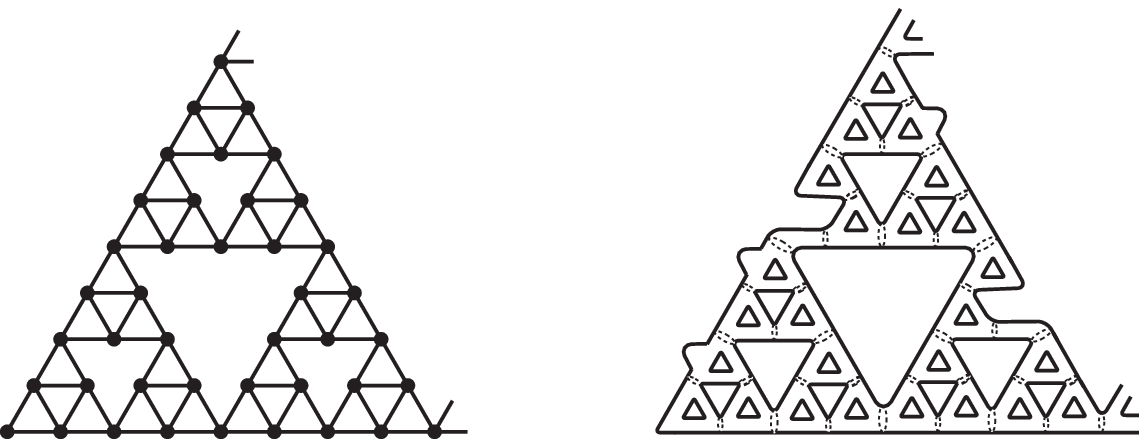}               

\ms
\cl{ Figure 1. Sierpinski gasket graph and fractal-like manifold}

\ms
The aim of this paper is twofold. First we  extend the results in 
[BB5, B2] to the  general framework of metric measure spaces. 
One of our main theorems gives necessary and sufficient
 conditions for  $\PHI(\Psi)$  to hold (Theorem 2.16). The
proofs given here are based closely on those in [BB5]. 
Secondly, we establish the stability of $\PHI(\Psi)$ under rough isometries
and under bounded perturbations, 
assuming some local regularity on the spaces (Theorem 2.21).
Let us look at two examples. 
The right side of Figure 1 is an example of a fractal-like manifold
$M$: it is a 
2-dimensional Riemannian manifold which is made from the Sierpinski gasket graph 
(left side of Figure 1) by replacing the edges by tubes of length $1$, 
and by gluing the tubes
together smoothly at the vertices. Also, one allows 
some small bumps, and some of the tubes to be removed.
More precisely, $M$ is a  
smooth Riemannian manifold which is roughly isometric 
to the Sierpinski gasket graph. 
Brownian motion moves on the surface of the tubes.
For our second example, Figure 2 is an example of a fractal tiling 
which is made by the following procedure. 
First, consider the triangular lattice on $\bR^2$ where each edge is of length $1$. 
We then fill each triangle with a Sierpinski gasket. This fractal tiling 
is roughly isometric to $\bR^2$. One can define a
Dirichlet form on the space by summing up Dirichlet forms on the gasket. 
As a consequence of our main theorems and the known heat kernel estimates on the 
Sierpinski gasket
([BP, Jo]), we have the following: the Dirichlet form for the manifold
depicted on the right side of Figure 1 satisfies 
$\PHI(\Psi)$ with $\Psi (s)=s^{2}\vee s^{\log 5/\log 2}$ while the 
Dirichlet form for the tiling shown in Figure 2 satisfies 
$\PHI(\Psi)$ with $\Psi (s)=s^{2}\wedge s^{\log 5/\log 2}$. 

\ms \qquad\qquad\qquad \epsfysize 1.8 in   \epsffile{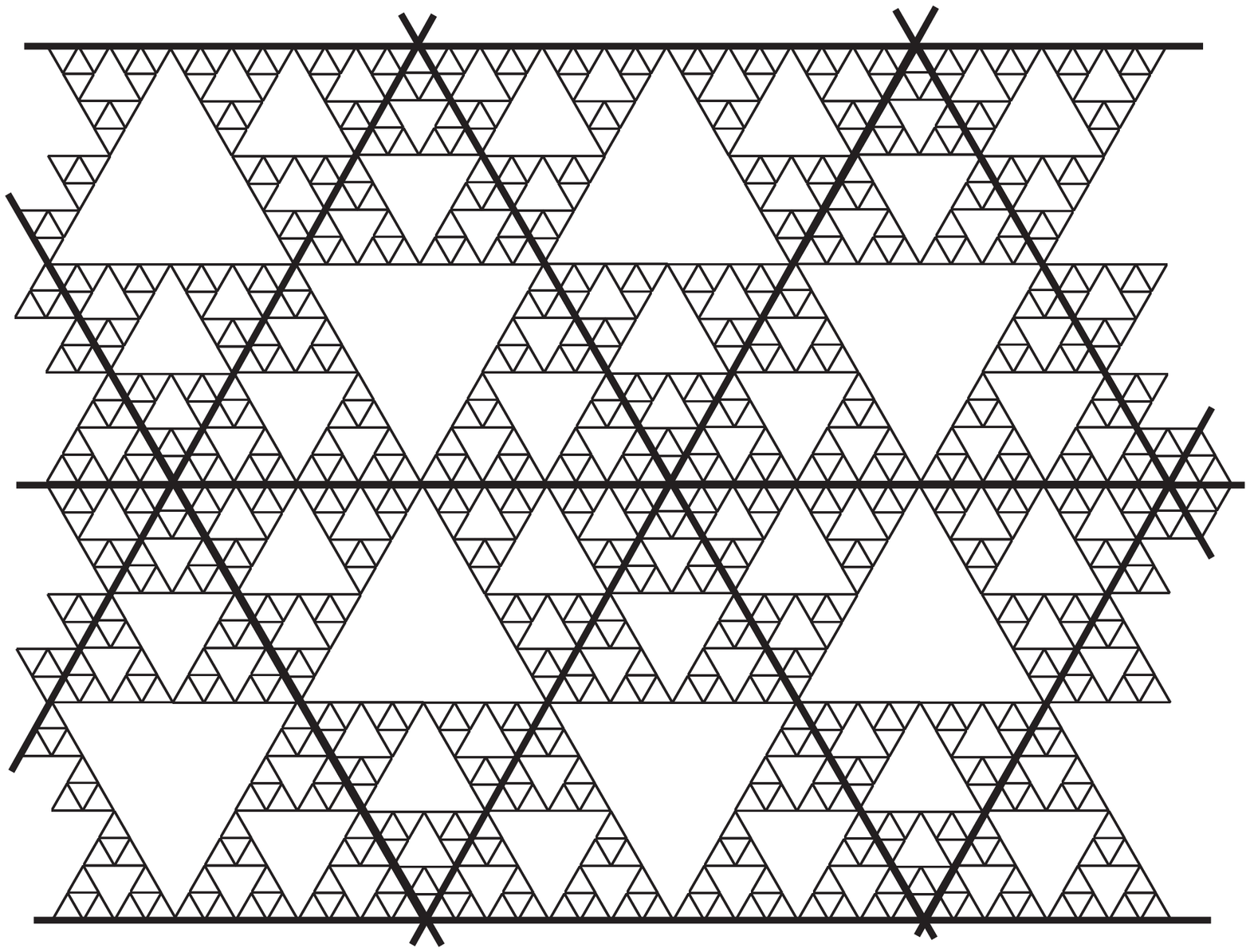}               

\ms
\cl{ Figure 2. Tiling of Sierpinski gaskets}

\ms
The organization of the paper is the following. In Section 2, we
present the framework and the main theorems. Sections 3 and 4 provide
the proof of Theorem 2.16. In Section 5, we give the proof of Theorem
2.21. We also give several applications. Our results, together with
Hino's results in [Hin], complete the proof of the singularity of the
energy measures of Dirichlet forms on higher dimensional Sierpinski
carpets with respect to the Hausdorff measures.  Another application
gives transition density estimates for the self-similar Brownian
motions defined on higher-dimensional Sierpinski carpets by [KZ].

We will use $c$, with or without subscripts, to denote strictly
positive finite constants whose values are insignificant and may
change from line to line.
\ms
\subsec{\bf 2. Framework and main theorems}

\ms 

Let $(X,d)$ be a connected locally compact 
complete separable metric space. 
We assume that the metric $d$ is geodesic: for each 
$x, y \in X$ there exists a (not necessarily unique) geodesic 
path $\gamma(x,y)$ such that for each $z\in \gamma(x,y)$, 
we have
$d(x,z)+d(z,y)=d(x,y)$.
Let $\mu$ be a Borel measure on $X$ such that $0<\mu(B)<\infty$ for
every ball $B$ in $X$. We write $B(x,r)=\{y: d(x,y)<r\}$,
and $V(x,r)=\mu(B(x,r))$.
Note that under the assumptions above, the closure of $B(x,r)$ is compact
for all $x\in X$ and $0<r<\infty$.  
For simplicity in what follows, we will also assume that
$X$ has infinite diameter, but similar results (with obvious 
modifications to the statements and the proofs) hold when 
the diameter of $X$ is finite. We will call such a space a {\sl 
metric measure space}, or a MM space.
`Metric measure space' is not a new term, and we mention 
[Gr2, GHL, Ha, Kas, Ke] as a sample of recent papers
that deal with analysis on such spaces.

Now let $(\sE,\sF)$ be a regular, 
strong local Dirichlet
form on $L^2(X,\mu)$: see [FOT] for details. We denote by $\Delta$ the 
corresponding self-adjoint operator; that is, we say $h$ is in the domain
of $\Delta$ and $\Delta h=f$ if $h\in \sF$ and $\sE(h,g)=-\int fg \, d\mu$
for every $g\in \sF$.
Let $\{P_t\}$ be the corresponding semigroup. 
$(\sE,\sF)$ is called conservative (or stochastically complete) if
$P_t1=1$ for all $t>0$. Throughout the paper, we assume that $(\sE,\sF)$ is 
conservative.

Since $\sE$ is regular, 
$\sE(f,g)$ can be written in terms of a signed
measure $\Gam({f,g})$. To be more precise, 
for $f\in \sF_b$ (the collection $\sF_b$  is the set of functions in $\sF$ that
are essentially bounded) $\Gam({f,f})$ is the unique
smooth Borel measure (called the energy measure) on $X$ satisfying 
$$\int_X{\tilde g}d\Gam({f,f})=2\sE (f,fg)-\sE (f^2,g),\qquad  g\in \sF_b,
$$
where ${\tilde g}$ is the quasi-continuous modification of $g\in \sF$. 
(Recall that $u: X\to \bR$ 
is called quasi-continuous if for any $\eps>0$, there exists an
open set $G\subset X$ such that $\hbox{Cap} (G)<\eps$ and $u|_{X\setminus G}$ is 
continuous. 
It is known that each $u\in \sF$ admits a quasi-continuous 
modification ${\tilde u}$ -- see [FOT], Theorem 2.1.3.) 
Throughout the paper, we will abuse notation and take the  quasi-continuous 
modification of $g\in \sF_b$ without writing ${\tilde g}$.
$\Gam({f,g})$ is defined by
$$\Gam({f,g})=\frac 12 (\Gam({f+g,f+g})-\Gam({f,f})
-\Gam({g,g})),\qquad  f,g\in \sF.$$ 
$\Gam({f,g})$ is also local, linear in 
$f$ and $g$, and satisfies 
the Leibniz and chain rules -- see [FOT], p. 115-116. 
That is, if $f_1, \ldots, f_m, g$, and $\vp(f_1, \ldots, f_m)$
are in $\sF_b$, and $\vp_i$ denotes the partial derivative
of $\vp$ in the $i^{th}$ direction, we have:
$$ d\Gam({fg,h})= fd\Gam(g,h)+gd\Gam(f,h), $$
$$ d\Gam(\vp(f_1,\dots, f_m),g)
 =\sum_{i=1}^m \vp_i(f_1,\dots, f_m) d\Gam(f_i,g). $$
We call $(X,d,\mu,\sE)$ a metric measure Dirichlet space, or a 
MMD space. 

\med {\bf Examples.}
\ni 1. If $M$ is a Riemannian manifold, we can take
$d$ to be the Riemannian metric and $\mu$ the 
Riemannian measure. The Dirichlet form $\sE$ is defined
by taking its core $\sC$ to be the $C^\infty$ functions
on $M$ with compact support, and defining
$$ \sE(f,f) =\int_M |\grad f|^2 d\mu, \q f\in \sC.$$
The domain $\sF$ of $\sE$ is then the completion of
$\sC$ with respect to the norm $||f||_2 + \sE(f,f)^{1/2}$, 
and $d\Gam(f,g) = \grad f \cdot\grad g\, d\mu$. 

\small 2. Cable system of a graph.
Given a weighted graph $(G,E,\nu)$ (see Definition 2.13 below)
we can define the {\sl cable system}
$G_C$ by replacing each edge of $G$ by a copy of $(0,1)$, joined
together in the obvious way at the vertices. 
For further details see [BB5] or [BPY].
Let $\mu$ be the measure on $G_C$
given by taking $d\mu(t)=\nu_{xy}\, dt$ for $t$ in the cable
connecting $x$ and $y$,
where $\nu_{xy}$ is the conductance of the edge connecting
$x$ and $y$; see [BB5]. 
One takes as the core $\sC$ the functions in $C(G_C)$ which
have compact support and are $C^1$ on each cable, and sets
$$ \sE(f,f) = \int_{G_C} |f'(t)|^2 d\mu(t). $$
One use of this construction is that the restriction to $G$
of a harmonic function $h$ on $G_C$ yields a harmonic function on $G$.

\small 3. Let $D$ be a domain in $\bR^d$ with a smooth boundary.
Then let $\sC=C_0^2(\ol D)$, $\mu$ be Lebesgue measure, and 
$$ \sE(f,f) = \half \int_D |\grad f|^2 d\mu. $$
The associated Markov process $Y$ is Brownian motion on $D$
with normal reflection on $\pd D$. 
For the extension of this construction to piecewise smooth domains
such as the pre-Sierpinski carpet, see [BB4].

\small 4. For fractal sets it is not as easy to describe
$\sE$. However, let  $F\subset \bR^d$ be a connected set,
and suppose that there exists a geodesic metric $d$ on $F$. 
Let $\mu$ be the Hausdorff 
$\al$-measure on $F$ (with respect to $d$) and suppose that
$$ c_1 r^\al \le \mu(B_d(x,r)) \le c_2 r^\al, \q x\in F, \, r>0. $$
Let 
$$ \eqalign{
 N_{\sigma, \infty}(f)
 &= \sup_{0<r\le 1} r^{-\al-2\sigma} \int_F\int_F 1_{B(y,r)}(x) 
  |f(x)-f(y)|^2 d\mu(x) d\mu(y), \cr
 \Lam^\sigma_{2, \infty}(f) 
 &= \{ u\in L^2(F,\mu):  N_{\sigma, \infty}(f)< \infty \}. \cr}$$
There exist many 
fractals satisfying the above with 
a Dirichlet form $\sE$ on $L^2(F,\mu)$ 
for which 
the domain $\sF$ of $\sE$ is given by $ \Lam^{\beta/2}_{2, \infty}$,
and $c_1 N_{\sigma, \infty}(f) \le \sE(f,f) \le c_2 N_{\sigma, \infty}(f)$;
see [Gr2, GHL, Kum2] etc. 

In the particular case of the Sierpinski gasket $F=F_{SG}$,
let $F_n$ be the set of vertices of triangles of side $2^{-n}$;
regard $F_n$ as a graph with $x\sim y$ if and only if $x$ and $y$
are in some triangle of side $2^{-n}$. 
Then for $f\in  \Lam^{\beta/2}_{2, \infty}$ 
with $\beta=\log 5/\log 2$, 
one has
$$ \sE(f,f) = c \lim_{n\to \infty} (5/3)^n \sum_{x\sim y} 
 (f(x)-f(y))^2 . $$
\ms
In many contexts (including examples 1--3 above) 
there is a natural metric $d_\sE$ associated with the
Dirichlet space $(\sE,\sF, L^2(M,\mu))$, defined 
as follows. Let $$\sL(\lam)= \{ f\in \sF:  d\Gam(f,f) \le \lam d\mu\},\eqno (2.1)$$
and set
$$ d_\sE(x,y) =\sup\{ f(x)-f(y) : f\in \sL(1) \}. $$
We cannot use this metric in this paper, since for many
true fractal spaces it is known 
(see [Kus, Hin]) that $\sL(1)$ contains only (quasi-everywhere) 
constant functions, so that $d_\sE(x,y) \equiv 0$. 
(One can also have, in other circumstances, $d_\sE(x,y) \equiv\infty$ -- see 
e.g. [BS]). 

\ms Let $\beta, {\betaloc}\ge 2$ and 
$$ \Psi(s)= \Psi_{\betaloc, \beta}(s) 
 = \cases {s^{\betaloc} & if $s\le 1$ \cr
 s^\beta & if $s> 1$. \cr }  
$$
$\Psi(s)$ will give the space/time scaling on the space $X$. 
Generalization of this time scaling factor may be
possible (see e.g., [T]), but we do not pursue it here. 

\ms We now give a number of conditions which $X$ may or may not 
satisfy.

\med {\bf Definition 2.1.} 
(a)  $X$ satisfies volume doubling $\VD$ 
if there exists a constant $c_1$ such that
$$ V(x,2R) \le c_1  V(x,R) \q \hbox{ for all } x\in X, \, R\ge 0. \eqno(\VD)$$
\ni
(b) $X$ satisfies the {\sl Poincar\'e inequality} $\PI(\Psi)$
if there exists a constant $c_{2}$ such that
for any ball $B=B(x,R) \subset X$ and 
$f\in \sF$, 
$$ \int_{B} (f(x)-\ol f_B)^2 d\mu(x)  \le c_{2} \Psi(R) 
  \int_B d\Gam ({f,f}) . \eqno (\PI(\Psi))$$
Here $\ol f_B = \mu(B)^{-1} \int_B f(x) d\mu(x)$.

\ni (c) 
We say a function $u$ is {\sl harmonic} on a domain $D$ if
$u\in \sF_{loc}$ and 
$\sE(u,g)=0$ for all $g\in \sF$ with support in $D$.
Here $u\in \sF_{loc}$ if and only if for any relatively compact open set
$G$, there exists a function $w\in \sF$ such that $u=w$ $\mu$-a.e. on $G$. 
See page 117 in [FOT] for the definition of $\sE(u,g)$ for $u\in \sF_{loc}$
when $(\sE,\sF)$ is a regular, strong local Dirichlet form.
Functions in $\sF$ are only defined up to quasi-everywhere equivalence;
we use a quasi-continuous modification of $u$.
$X$ satisfies the {\sl elliptic Harnack inequality} $\EHI$ if
there exists a constant $c_3$ such that, 
for any ball $B(x,R)$, whenever $u$
is a non-negative harmonic function on $B(x,R)$ then
there is a quasi-continuous modification $\tilde u$ of $u$ that 
satisfies 
$$ \sup_{B(x,R/2)} {\tilde u} \le c_3 \inf_{B(x,R/2)} {\tilde u}. \eqno(\EHI)$$ 
Note that by a standard argument (see, e.g., [M1], p.~587) $\EHI$ implies that 
$\tilde u$ is H\"older continuous.

\ni (d) 
Let  $Q=Q(x_0,T,R)=(0, 4T)\times B(x_0, 2R)=:I\times B_{2R}$. 
Let $u(t,x): Q \to \bR$. 
\item{$\bullet$}  We define
$u_t=\frac{\del u}{\del t}\in L^2(dt\times \mu)$ as the derivative in
the Schwartz' distribution sense. That is, we
define $u_t$ to be the function $f$  in $L^2(dt\times \mu)$ so that for any
function $g: Q \to \bR$ such that $g(x,\cdot) \in C^\infty_K(0,4T)$
for each $x\in	B(x_0, 2R)$ and
$g_t = \frac{\pd g}{\pd t}\in L^2(dt\times \mu)$,
then
$$ \int_Q ( f(x,t)g(x,t) + u(x,t)g_t(x,t))\,dt\,\mu(dx)=0. $$
\item{$\bullet$} Let $H(I\to \sF^*)$ be the space of functions 
$u\in L^2(I\to\sF^*)$ with the distributional time derivative 
$u_t\in L^2(I\to\sF^*)$ equipped with the norm 
$$\Big(\int_I\|u(t,\cdot)\|_{\sF^*}^2+\|u_t(t,\cdot)\|_{\sF^*}^2dt\Big)^{1/2}.$$
Here we identify $L^2(X,\mu)$ with its own dual and denote the dual of $\sF$ by 
$\sF^*$. So, $\sF\subset L^2(X,\mu)\subset \sF^*$ with continuous and dense embeddings. 

Let $\sF (I\times X)=L^2(I\to\sF)\cap H(I\to \sF^*)$ be a Hilbert space with norm 
$$\|u\|_{\sF(I\times X)}=
\Big(\int_I\|u(t,\cdot)\|_{\sF}^2+\|u_t(t,\cdot)\|_{\sF^*}^2dt\Big)^{1/2}.$$
\item{$\bullet$} We define $\sF_{loc}(Q)$ to be the set of 
$dt\otimes d\mu$-measurable functions on $Q$ such that for every relatively
compact open set  
$D\subset\subset B_{2R}$ and every open interval $I'\subset\subset I$, 
there exists 
a function $u'\in \sF (I\times X)$ with $u=u'$ on $I'\times D$. We define
$$\sF_c(Q):=\{u\in \sF (I\times X): u(t,\cdot) \hbox{ has compact support
in $B_{2R}$ for a.e. } t\in I\}.$$

We say a function
$u(t,x): Q \to \bR$ is a solution of the heat equation in $Q$ if
$u\in \sF_{loc}(Q)$ and 
$$\int_J \Big[\int f(t,x) u_t(t,x) \mu(dx) 
+\sE(f(t,\cdot), u(t,\cdot))\Big] dt=0,~~~
\forall J\subset\subset I,~\forall f\in \sF_c(Q). \eqno (2.2) $$ 
$X$ satisfies the {\sl parabolic Harnack inequality} $\PHI(\Psi)$,
if there  exists a constant $c_{4}$ such that the following holds.
Let $x_0\in X$, $R>0$, $T=\Psi(R)$, and 
$u=u(t,x)$ be a non-negative solution of the heat equation
in $Q(x_0,T,R)$.
Write  $Q_- = (T,2T)\times B(x_0,R)$ and $Q_+ = (3T,4T)\times B(x_0,R)$;
then there exists $\tilde u=\tilde u(t,x)$  
such that $\tilde u(t, \cdot)$ is a quasi-continuous modification of
$u(t,\cdot)$ for each $t$ and
$$ \sup_{Q_-} {\tilde u} \le c_4 \inf_{Q_+} {\tilde u}. \eqno(\PHI(\Psi))$$
Given this PHI, a standard oscillation argument implies that
$\tilde u$ is jointly continuous.

\small
\ni{\bf Remark 2.2.} In the case of general MMD spaces we can only
define harmonic functions up to quasi-everywhere equivalence.
This is why we needed to be careful in our definitions of EHI
and $\PHI(\Psi)$.
\ms
It is easy
to deduce from $\VD$ that there exist $c_5, \al<\infty$ such that if $x, y
\in X$ and $0<r <R$ then
$$ \frac{V(x,R)}{V(y,r)} 
\le c_5 \Bigl( \frac{d(x,y)+R}{r} \Bigr)^{\al}. \eqno(2.3)$$

The following covering lemma holds; cf.~[BB5]. 

\proclaim Lemma 2.3. Assume that $X$ satisfies $\VD$. 
For $x_0\in X$ and $0<s\le R\le \infty$, there exists a cover 
of $B(x_0,R)$ by balls $B(x_i,s)$ with $x_i\in B(x_0,R)$ such that no point in 
$X$ is in more than $L_0$ of the $B(x_i,2s)$. Here $L_0$ depends only on 
$s$, $R$ and the constant $c_1$ in $\VD$.

\med {\bf Definition 2.4.} Let $A$, $B$ be disjoint subsets
of $X$. We define the effective resistance
$\Reff(A,B)$ by
$$ \Reff(A,B)^{-1} = \inf \Bigl\{  \int_X d\Gam ({f,f}) : 
 f=0 \hbox { on $A$ and $f=1$ on B}, \, f\in \sF \Bigr\}.\eqno(2.4)$$
$X$ satisfies the condition $\RES(\Psi)$ if there exist
constants $c_1, c_2$ such that for any $x_0\in X$,
$R\ge 0$,
$$ c_1 \frac{\Psi(R)}{V(x_0,R)}\le \Reff(B(x_0,R), B(x_0,2R)^c) \le 
c_2 \frac{\Psi(R)}{V(x_0,R)}. \eqno(\RES(\Psi))$$

\med {\bf Definition 2.5.}
$X$ satisfies $\CS(\Psi)$ 
if there exist  $\th \in (0,1]$ and constants $c_1, c_2$ such that
the following holds.
For every $x_0 \in X$, $R>0$ there exists a cut-off function
$\vp (= \vp_{x_0,R})$ with the properties:
\nl (a) $\vp(x)\ge 1$ for $x\in B(x_0,  R/2)$.
\nl (b) $\vp(x)=0$ for $x\in B(x_0, R)^c$.
\nl (c) $|\vp(x)-\vp(y)| \le c_1 (d(x,y)/R)^\th$ for all $x$,$y$.
\nl (d) For any ball $B(x,s)$
with $0<s\le R$ and 
$f\in \sF$,  
$$ \int_{B(x,s)} f^2 d\Gam(\vp,\vp) \le 
 c_2 (s/R)^{2\th} \Big( \int_{B(x,2s)} d\Gam(f,f) 
+ \Psi(s)^{-1} \int_{B(x,2s)} f^2  d\mu\Big).\eqno (2.5)$$

\med {\bf Remarks 2.6.} 
\nl  1. We call (2.5) a weighted Sobolev inequality. It is clear  
that to prove (2.5) it is enough to consider nonnegative $f$.
\nl 2. Suppose $\CS(\Psi)$ holds for
$X$, but with (a) above replaced by 
$$  \vp(x)\ge 1 \hbox{ for } x\in B(x_0, \delta R), \eqno(2.6)$$
for some $\delta <\half$.
Then an easy covering argument (using $\VD$) gives $\CS(\Psi)$
with $\delta=\half$.
\small 3. Let $\lam>1$. 
Suppose that $\CS(\Psi)$ holds, except that
instead of (2.5) we have        
$$ \int_{B(x,s)} f^2 d\Gam(\vp,\vp) \le 
 c_2 (s/R)^{2\th} \Big( \int_{B(x,\lam s)} d\Gam(f,f) 
+ \Psi(s)^{-1} \int_{B(x,\lam s)} f^2  d\mu\Big). \eqno (2.7)$$
Then once again it is easy to obtain  $\CS(\Psi)$ with $\lam=2$
by a covering argument.
\small 4. Any operation on the cut-off function $\vp$ which reduces
$d\Gam(\vp,\vp)$ while keeping properties 
(a), (b) and (c) of Definition 2.5  will generate a new cut-off function
which still satisfies (2.5).
We can therefore assume that any cut-off function $\vp$ 
satisfies the following:
\i{(a)} $0\le \vp \le 1$.
\i{(b)} For each $t\in (0,1)$ the set $\{x: \vp(x) >t\}$ is connected
 and contains $B(x_0,R/2)$.
\i{(c)} Each connected component $A$ of $\{x: \vp(x) < t\}$
intersects $B(x_0,R)^c$.
\small 5. Note that if $\CS(\Psi)$ holds for $\Psi= \Psi_{\betaloc, \beta}$, 
then $\CS(\Psi_{\betaloc',\beta'})$ holds if $\beta'\ge \beta$
and $\betaloc' \le \betaloc$.

\ms We have the following by the same arguments as in Lemma 5.1 of [BB5].

\proclaim Lemma 2.7. 
Let $X$ satisfy $\VD$, $\PI(\Psi)$ and $\CS(\Psi)$.
Then $X$ satisfies $\RES(\Psi)$.

\ms In the definitions which follow, recall that we have a 
`crossover' from $t=r^{\betaloc}$ scaling to $t=r^\beta$ scaling
at $r=1$. For $(t,r)\in (0,\infty)\times [0,\infty)$
we consider the two regions:
$$ \Lam_1 = \{ (t,r): t\le 1 \vee r \},~~~~ 
\Lam_2 = \{ (t,r): t\ge 1 \vee r \}. $$

Let 
$$ h_\beta(r,t) = \exp\Big(- \Big(\frac{r^\beta}{t}\Big)^{1/(\beta-1)}\Big). $$

\med {\bf Definition 2.8.} Let $X$ be a MMD space. We say $X$
satisfies $\HK(\Psi)$ if the heat kernel $p_t(x,y)$ on $X$ exists and
satisfies
$$ \frac{c_1 h_{\betaloc}( c_2 d(x,y),t) }{V(x,t^{1/{\betaloc}})} 
\le p_t(x,y) \le 
 \frac{c_3 h_{\betaloc}( c_4 d(x,y),t) }{V(x,t^{1/{\betaloc}})},\eqno(2.8)$$
for $x,y \in X$ and $t \in (0,\infty)$ with $(t,d(x,y))\in \Lam_1$, and 
$$ \frac{c_1 h_{\beta}( c_2 d(x,y),t) }{V(x,t^{1/{\beta}})} 
\le p_t(x,y) \le 
\frac{c_3 h_{\beta}( c_4 d(x,y),t) }{V(x,t^{1/{\beta}})}, \eqno(2.9)$$
for $x,y \in X$ and $t \in (0,\infty)$ with $(t,d(x,y))\in \Lam_2$. 

\med {\bf Remark 2.9.} 
\nl To understand why the crossover takes the form it does, it is useful
to consider the contribution to $p_t(x,y)$ from various types
of paths in $X$. Let $r=d(x,y)$. 
First, if $0<t\le 1$ and $r<1$ then the behaviour is essentially local.

If $r\geq t$ then we are in the `large deviations' regime:
the main contribution to $p_t(x,y)$ is from those
paths of the Markov process $Y$ which 
are within a distance $O(t/r)$ of a geodesic from $x$ to $y$. 
So, once the length of the geodesic is given, only the local 
structure of $X$ plays a role.
Note that in this case the term in the exponential is smaller
than $e^{-c t}$, so that the volume term $V(x,t^{1/{\betaloc}})^{-1}$ could
be absorbed into the exponential with a suitable modification
of the constants $c_2$ and $c_4$.

Finally, if $t>1$ and $r<t$, then the paths which contribute to
$p_t(x,y)$ fill out a much larger part of $X$: those which lie in 
$B(x,t^{1/\beta})$ if $r<t^{1/\beta}$, and 
those which are within a distance $O(t/r^{\beta-1})$ of a geodesic 
from $x$ to $y$ in the case when $t^{1/\beta} \le r \le t$.

\ms We will also want to discuss local versions of these conditions.
Since these will only depend on
$\Psi(s)$ for $s\in [0,1]$, they are independent of the 
the parameter $\beta$.
We say  $X$ satisfies $\VD_{\loc}$ if 
$\VD$  
holds for $x\in X$, $0<R\le 1$.
Similarly we define $\PI(\bar\beta)_\loc$, $\EHI_\loc$, $\CS(\bar\beta)_\loc$ and 
$\PHI(\bar\beta)_\loc$ by requiring the conditions only for $0<R\le 1$.
For $\HK(\bar\beta)_\loc$ we require the bounds only for
$t\in (0,1)$ -- so only (2.8) is involved.
The value 1 here is just for simplicity: each of the local
conditions implies an analogous local condition for $0<R\le R_0$
for any (fixed) $R_0 >1$ -- see Section 2 of [HSC].

Finally, we introduce two local notions which do not include any scaling order.

\med {\bf Definition 2.10.} (a) We call $\vp$ a {\sl cut-off function} for 
$A_1 \subset A_2$ if $\vp =1 $ on $A_1$ and is zero on $A_2^c$.

\ni (b) We say  $X$ satisfies $\PI_\loc$ if 
for each $c_1>0$, there exists $c_2>0$ such that 
$$ \int_{B} (f(x)-\ol f_B)^2 d\mu(x)  \le c_{2} 
\int_B d\Gam ({f,f}) \eqno(2.10)$$
for any ball $B=B(x,c_1) \subset X$ and $f\in \sF$.  

\ni (c) We say $X$ satisfies $\CC_\loc$ 
if for every $x_0 \in X$, there exists a cut-off function
$\vp (= \vp_{x_0})$ for $B(x_0,  1/2)\subset B(x_0,  1)$ such that 
$$ \int_{B(x_0,1)} d\Gam(\vp,\vp) \le c_3 V(x_0,1),$$
where $c_3>0$ is independent of $x_0$ and $\vp$.

\ni $\CC$ stands for `controlled cut-off' functions.
Clearly $\PI(\betaloc)_\loc$ for any $\betaloc\ge 2$ implies 
$\PI_\loc$ and  $\CS(\betaloc)_\loc$ for any $\betaloc>0$ implies 
$\CC_\loc$.

\smallskip We next give some sufficient condition for $\CS(2)_\loc$.

We say $(X,d,\mu, \sE)$ has {\sl regular
local cut-off functions} if there exists a constant $C_1$ 
such that for all $x\in X$, $R\in (0,1)$ there exists a 
cut-off function $\vp:X \to \bR$ for 
$B(x,R/2)\subset B(x,R)$ such that $R \vp \in \sL(C_1)$, 
where $\sL(C_1)$ is defined in (2.1). 

\proclaim Lemma 2.11.
If the MMD space $(X,d,\mu, \sE)$ has regular
local cut-off functions such that 
(c) of Definition 2.5 holds, 
then it satisfies $\CS(2)_\loc$.

\proof Let $\vp$ be a regular local
cut-off function  for 
$B(x,R/2)\subset B(x,R)$ such that $R \vp \in \sL(C_1)$.
Then $d\Gam(\vp,\vp) \le C_1^2 R^{-2} d\mu$, so
$\int_{B(x,s)} f^2 d\Gam(\vp,\vp) \le 
 C_1^2 R^{-2} \int_{B(x,s)} f^2 d\mu$. 
Thus (2.5) holds. \qed

If $d_\sE$ is a true metric, then  $(X,d_\sE,\mu, \sE)$ 
has  regular local cut-off functions since one can take $C_1=2$ and 
$\vp(y)=(1 \wedge (2- 2R^{-1} d_\sE(x_0,y))\vee 0$.
Thus Examples 1--3 above have regular local cut-off functions. 

\med {\bf Definition 2.12.} 
$X$ satisfies the condition $\EE(\Psi)$ if for any $x_0\in X$,
$R\ge 0$,
$$c_1 \Psi(R)\le \E^{x_0}[\tau_{B(x_0,R)}] \le c_2 \Psi(R), \eqno(\EE(\Psi))$$
where $\tau_{A}=\inf\{t\ge 0: Y_t\notin A\}$, $Y_t$ is the strong
Markov process associated to the Dirichlet form $(\sE, \sF)$,
and $\E^{x_0}$ denotes the expectation starting from the point $x_0$.
$\EE(\Psi)$ describes the `walk dimension' of the associated Markov process.

We also will need to consider weighted graphs. 

\med{\bf Definition 2.13.} 
Let $(G,E)$ be an infinite locally finite connected graph. 
We write $x\sim y$ if $(x,y)\in E$, i.e., there is an edge connecting
$x$ and $y$.
Define edge weights (conductances) 
$\nu_{xy}=\nu_{yx}\ge 0$, $x,y\in G$, and assume that
$\nu$ is adapted to the graph structure by requiring that
$\nu_{xy}>0$ if and only if $x\sim y$.
Let $\nu_x = \sum_y \nu_{xy}$, and define a measure $\nu$
on $G$ by $\nu(A)=\sum_{x\in A} \nu_x$. 
We call $(G,\nu)$ a {\sl weighted graph}.  

We write $d(x,y)$ for the graph distance, and define the balls
$$ B_G(x,r)=\{y: d(x,y) < r\}.$$ 
Given $A\subset G$ write 
$\partial A = \{ y\in A^c: d(x,y)=1$ for some $x \in A\}$ for the exterior
boundary of $A$, and let $\ol A=A \cup \partial A$.

\med {\bf Definition 2.14.} A weighted graph 
$(G,\nu)$ has {\sl controlled weights}
if there exists $p_0>0$ such that for all $x,y\in G$ 
$$ \frac{\nu_{xy}}{\nu_x} \ge p_0, \qq x\sim y . $$
This was called the $p_0$-condition in [GT2].

\ms The Laplacian is defined on  $(G,\nu)$ by
$$ \Delta f(x) = \frac{1}{\nu_x} \sum_y \nu_{xy} (f(y)-f(x)). $$
We also define a Dirichlet form $(\sE,\sF)$ by taking
$\sF=L^2(G,\nu)$, and
$$ \sE(f,g) = \half \sum_x \sum_y (f(x)-f(y))(g(x)-g(y))\nu_{xy},
\q f,g \in \sF. $$
If $f \in \sF$  we define the measure $\Gam_G(f,f)$ on $G$
by setting
$$ \Gam_G(f,f)(x) = \sum_{y\sim x}  (f(x)-f(y))^2 \nu_{xy}. $$

The conditions $\VD$, $\EHI$ and  $\PHI(\Psi)$ for graphs are defined in exactly the
same way as for manifolds; see [BB5]. 
The definitions of  $\PI(\Psi)$ and $\RES(\Psi)$ are also the
same. For the bound $\HK(\Psi)$ we only require (2.9).
The condition $\CS(\Psi)$ is also the same; the weighted Sobolev
inequality (2.5) takes the form
$$ \eqalignno{ 
 \sum_{x\in B_G(x_1,s)} &f(x)^2 \Gam_G(\vp,\vp)(x)  &(2.11) \cr
 &\le   c_2 \bigl( \frac{s}{R} \bigr)^{2\th} \Big( 
 \sum_{x\in B_G(x_1,2s)}  \Gam_G(f,f)(x)
   + \Psi(s)^{-1} \sum_{x\in B_G(x_1,2s)} \nu_x f(x)^2 \Big). }$$

\ni It is easy to check that $\PI_\loc$ and $\CC_\loc$ hold for any weighted 
graph with controlled weights. In fact, $\PI(\betaloc)_\loc$ and $\CS(\betaloc)_\loc$
hold for any choice of $\betaloc\ge 2$ on such graphs, since it is irrelevant to treat
$R<1$ for graphs.

\bs
We summarize the conditions we have introduced: 
\medskip

\halign{\qquad \qquad#\hfil \qquad \qquad &#\hfil \qquad \qquad &#\hfil \cr
VD&Volume doubling\cr
$PI(\Psi)$ & Poincar\'e inequality \cr
EHI&Elliptic Harnack inequality\cr
$\PHI(\Psi)$ & Parabolic Harnack inequality \cr
$RES(\Psi)$&Resistance exponent \cr
CS($\Psi$)&Cut-off Sobolev inequality\cr
$\CC$&Controlled cut-off functions\cr
HK($\Psi$)& Heat kernel estimates\cr
$\EE(\Psi)$&Walk dimension\cr
}
\bs

\smallskip
We will need the following:

\proclaim Theorem 2.15. 
(See [BBK], [HSC] Theorem 5.3, [GT3]).
The following are equivalent:
\nl (a) $X$ satisfies $\PHI(\Psi)$.
\nl (b) $X$ satisfies $\HK(\Psi)$.
\nl (c) $X$ satisfies $\VD$, $\EHI$ and $\RES(\Psi)$. 
\nl (d) $X$ satisfies $\VD$, $\EHI$ and $\EE(\Psi)$. 

\proof The equivalence of (a) and (b) is given in [BBK]. 
(See also [HSC] for the case where solutions to the heat equation are
sufficient regular.)  That these are equivalent to (c) and (d) is
proved in [GT3].  (See [GT2] for the graph case.)  \qed 

\ms 
The first of our main theorems is the following. (The graph case was
proved in [BB5]).

\proclaim Theorem 2.16. 
Suppose that $X$ is either an infinite connected weighted graph
with controlled weights, or a MMD  space. 
The following are equivalent:
\nl (a) $X$ satisfies $\VD$, $\PI(\Psi)$ and $\CS(\Psi)$.
\nl (b) $X$ satisfies $\PHI(\Psi)$.

Our second topic is the stability of $\PHI(\Psi)$. We will actually
discuss two kinds of stability. 

\med {\bf Definition 2.17.} A property $P$ is {\sl stable under 
bounded perturbation} if whenever $P$ holds for $(\sE^{(1)},\sF)$, then 
it holds for $(\sE^{(2)},\sF)$, provided 
$$c_1\sE^{(1)}(f,f)\le \sE^{(2)}(f,f)\le c_2\sE^{(1)}(f,f),~~~~\hbox{for all }
f\in \sF.\eqno(2.12)$$

The following result is due to Le Jan ([LJ], Proposition 1.5.5(b)). 
A simple proof is given in [Mos] p. 389.

\proclaim Lemma 2.18. Let $X$ be a MMD space. 
Suppose $(\sE^{(i)}, \sF), i=1,2$ are 
strong local regular Dirichlet forms  
that satisfy (2.12). Then the energy measures $\Gam^{(i)}$ satisfy
$$ c_1 d\Gam^{(1)}(f,f)\le d\Gam^{(2)}(f,f)\le c_2 d\Gam^{(1)}(f,f),
~~~~\hbox{for all }
f\in \sF.\eqno(2.13)$$

\smallskip
It is immediate from Lemma 2.18 that the conditions $PI(\Psi)$ and 
$CS(\Psi)$ are stable under bounded perturbations. So we deduce:

\proclaim Theorem 2.19. Let $X$ be a MMD space. 
Then $\PHI(\Psi)$ and $\HK(\Psi)$ are stable under 
bounded perturbations.

The second kind of stability is stability under rough isometries. 

\med {\bf Definition 2.20.} For each $i=1,2$, let $(X_i,d_i,\mu_i)$ be either 
a metric measure space or a weighted graph. A map
$\vp: X_1 \to X_2$ is a {\sl rough isometry} if there exist constants $c_1>0$ and 
$c_2, c_3>1$ such that
$$ X_2=\bigcup_{x\in X_1} B_{d_2}(\vp(x),c_1), $$
$$ c_2^{-1} (d_1(x,y)-c_1) \le d_2(\vp(x),\vp(y))\le  c_2 (d_1(x,y)+c_1),  $$
and
$$ c_3^{-1} \mu_{1}(B_{d_1}(x,c_1)) \le  \mu_2(B_{d_2}(\vp(x),c_1)) 
 \le  c_3 \mu_{1}(B_{d_1}(x,c_1)). $$
If there exists a rough isometry between two spaces they are
said to be {\sl roughly isometric}. (One can check this is an
equivalence relation.) 

\ms 
This concept was introduced by Kanai in [Kan1, Kan2].  A rough
isometry between $X_1$ and $X_2$ means that the global structure of
the two spaces is the same. However, to have stability of Harnack
inequalities, we also require some control over the local structure.
In the case of graphs it is enough to have controlled weights, but for
metric measure spaces more regularity is needed.  (In
[Kan1, Kan2] this local
control was obtained by geometrical assumptions on the manifolds).

\bs
The following theorem concerns the stability of $\PHI(\Psi)$  
under rough isometries. 

\proclaim Theorem 2.21. Let $X_i$ be either a MMD space satisfying 
$\VD_\loc$ and $\PI_\loc$ or a graph with controlled weights, 
and suppose there exists a rough isometry $\vp: X_1 \to X_2$. 
Let $\Psi_i(s)=s^{\bar\beta_i}1_{\{s\le 1\}}+s^\beta 1_{\{s\ge 1\}}$. 
\nl (a) Suppose that  $X_2$ satisfies $\PI(\bar\beta_2)_\loc$.
If $X_1$ satisfies $\VD$, 
$\CC_\loc$ and $\PI(\Psi_1)$ then 
$X_2$ satisfies $\VD$ and $\PI(\Psi_2)$.
\nl (b) Suppose that $X_2$ satisfies $\CS(\bar\beta_2)_\loc$.
If $X_1$ satisfies $\VD$ and $\CS(\Psi_1)$ then 
$X_2$ satisfies $\VD$ and $\CS(\Psi_2)$.

\smallskip By this theorem together with Theorem 2.16, we see that 
$\PHI(\Psi)$ is stable under rough isometries, given suitable local
regularity of the two spaces.

\subsec{3. Construction of cut-off functions}

\bs In this section we will prove 
Theorem 2.16 (b) $\Rightarrow$ (a). 
If $X$ satisfies PHI($\Psi$), then (VD) and HK($\Psi$) hold by Theorem 2.15.

\proclaim Lemma 3.1. PI($\Psi$) holds. 

\proof Let  $\ol Y$ be the process $Y$ 'reflected at the boundary' of 
the ball $B=B(x_0,R)\subset X$: see  [Ch] for
this construction in a general context. 
Write $\ol p_t(x,y)$ for the transition density of $\ol Y$,
and $p^0_t(x,y)$ for that of the process $Y$ killed on exiting $B$.

The argument of [SC1] gives $\PI(\Psi)$ once we have the estimate
$$  \ol p_t(x,y) \ge p^0_t(x,y) \ge c V(x_0,R)^{-1}, 
 \q x,y \in B(x_0,R/2), \, \half \Psi(R)\le t\le \Psi(R).  $$
This is proved in [BBK], as a key step in the derivation 
of HK($\Psi$) from PHI($\Psi$). \qed

\ms 
We now prove that $\PHI(\Psi)\Rightarrow \CS(\Psi)$, and
follow the arguments in [BB5] and [B2].  
The main difference from [BB5] is that a strong transience condition
(called (FVG)) was needed there in the initial arguments.  We will
also slightly simplify the arguments in [B2].

Let $D=B(x_0,R-\eps)$ where $\eps< R/10$, and $\lam>0$.
Let $Y$ be the process associated with the Dirichlet form
$(\sE, \sF)$.
Let $G_\lam^D$ be the resolvent associated with the process $Y$
killed on exiting $D$; that is,
$$ G_\lam^D f(x) = \bE^x \int_0^{\tau_D} e^{-\lam t} f(Y_t) dt, $$
for bounded measurable $f$, where $ \tau_D =\inf\{ t: Y_t \in X-D\}$. 
Let $p^D_t(\cdot,\cdot)$ be the heat kernel of $Y$ killed on exiting $D$.
Then the Green kernel of $G_\lam^D$ is given by
$$ g_\lam^D(x,y) =\int_0^\infty e^{-\lam t} p^D_t(x,y) dt.$$
We use the Green kernel to build a cut-off function $\vp$. 

\proclaim Lemma 3.2. Let $x_0\in X$. Then there exists
$\delta>0$ such that if $\lam= c_0  \Psi (R)^{-1}$  
$$ g^D_\lam(x_0,y) \le C_1 \frac{\Psi(R)}{V(x_0,R)}, \qq y \in B(x_0,
\delta R)^c, \eqno(3.1)$$
$$ g^D_\lam(x_0,y)  \ge C_2 \frac{\Psi(R)}{V(x_0,R)}, 
 \qq y \in B(x_0,\delta R).  \eqno(3.2)$$

\proof This follows easily from $\HK(\Psi)$ by integration. \qed

\proclaim Lemma 3.3. Let $x_0$ and $R$ be as above, and let 
$x,y\in B(x_0,\delta R)^c$. Then there exists $\th>0$ such that 
$$ |g^D_\lam(x_0,x)-g^D_\lam(x_0,y)| 
 \le c_1\Big(\frac{d(x,y)}{R}\Big)^{\th} \sup_{B(x_0,\delta R)^c} 
g^D_\lam(x_0,.). \eqno(3.3)$$

\proof The H\"older continuity of $p^D_t$ follows from $\PHI(\Psi)$
by a standard argument; see [M2].
Integrating we obtain (3.3).\qed

Fix $x_0\in X$ and let $B'=B(x_0,\delta R)$, $B=B(x_0,R)$, 
$D=B(x_0,R-\eps)$ where $\eps< R/10$.  Let $\lambda =c_0\Psi(R)^{-1}$ and 
define 
$$\vp (x)=1\wedge (c\Psi (R)^{-1}G^D_{\lambda}1_{B'}(x)),$$
where $c$ is chosen so that $\vp(x)=1$ on $x\in B'$. 
Using Lemmas 3.2 and 3.3, it is easy to check that 
$\vp$ is a cut-off function for $B'\subset B$ that satisfies 
Definition 2.5 (a)--(c). 
To complete the proof of $\CS(\Psi)$, we need to establish (2.5).

\proclaim Proposition 3.4. Let $x_1\in X$ and $f\in \sF$. 
Let $\delta$ be defined by Lemma 3.2 and let
$I=B(x_1,\delta s)$ with $0<s\le R$ and $I^*=B(x_1,s)$.  
There exist $c_1,c_2>0$ such that
for all $f\in \sF$,  
$$\int_{I} f^2 d\Gam(\vp,\vp) \le c_1 (s/R)^{2\th} 
 \Big(  \int_{I^*} d\Gam(f,f) 
 + c_2 \Psi(s)^{-1} \int_{I^*} f^2d\mu \Big). \eqno(3.4) $$

\proof {\bf Case 1.}  We first consider the case where $s=R$
and $x_1=x_0$. 
Let 
$$ \sF_D = \{ f\in \sF: \wt f=0 \hbox { q.e. on } X-D \}. $$
Set
$$\sE_\lam(f,g)=\sE(f,g)+\lam \int fg\, d\mu.$$
Let $v= G_\lam^D 1_{B'}$. Note that 
$$ v(x) \le \int_{B'}g^D(x,y)d\mu (y)\le \E^x[\tau_D]\le c \Psi(R), \qq x \in D,
\eqno(3.5)$$ by Theorem 2.15.
By [FOT] Theorem 4.4.1, $v\in  \sF_D$ and is quasi-continuous.
Further, since $Y$ is continuous, $v=0$ on $\ol D^c$.
Let $f\in \sF$. Then
$$ \int_B f^2 d\Gam(v,v) \le \int_X f^2 d\Gam(v,v) = 
 \int_X  d\Gam(f^2 v,v)   -\int_X 2fv  d\Gam(f,v). $$
Since $v \in \sF_D$ we have $f^2 v \in \sF_D$, so by
[FOT], Theorem 4.4.1,
$$  \int_X  d\Gam(f^2 v,v) = \sE(f^2 v,G^D_\lam 1_{B'}) 
 \le  \sE_\lam(f^2 v,G^D_\lam 1_{B'}) = \int_X f^2 v 1_{B'} d\mu
 \le c \Psi(R) \int_{B'} f^2 d\mu, $$
where we used (3.5) in the last inequality. 
Using Cauchy-Schwarz and (3.5), we obtain
$$ \eqalign{
 \Bigl|\int_X 2fv  d\Gam(f,v)\Bigr| &\le 
 c \Big(\int_X v^2 d\Gam(f,f) \Big)^{1/2} \Big(\int_X f^2 d\Gam(v,v) \Big)^{1/2} \cr
 &\le c \Psi(R) \Big(\int_B d\Gam(f,f) \Big)^{1/2} 
\Big(\int_X f^2 d\Gam(v,v) \Big)^{1/2}.\cr}$$

So, writing $H= \int_X f^2 d\Gam(v,v)$,  $J=\int_B d\Gam(f,f)$,
$K= \int_B f^2 d\mu$, we have
$$ H \le c \Psi(R) K + c \Psi(R) J^{1/2} H^{1/2}, $$
from which it follows that $H \le c \Psi(R) K + c\Psi(R)^2 J$.
{}From this, (3.4) with $s=R$ follows easily. 

\nl {\bf Case 2.} Define 
$$ Q(b)= Q(x_0,b) = \{ y: g^D_\lam (x_0,y) > b \}. $$
and let $$h=C_2\Psi(R)/(2V(x_0,R)),$$ 
where $C_2$ is as in Lemma 3.2. 
Note that by Lemma 3.2 and the fact  $g^D_\lam (x_0,y)=0$ 
for $y\notin D$, 
$$B(x_0,\delta R)\subset Q(2h)\subset Q(h)\subset B(x_0,R).$$
In Case 2, we will consider the situation that either 
$$ I^* \subset Q(2h) \eqno (3.6) $$
or 
$$  I^* \cap B(x_0,\delta R/2) =\emptyset  \eqno (3.7) $$
hold. Since $\vp\equiv 1$ on $Q(2h)$, (3.4) is clear
if (3.6) holds. 
Thus, we consider when (3.7) holds. Let 
$\psi_s (x)=1\wedge (c\Psi (s)^{-1}G^{B(x_0,s-\epsilon)}_{\lambda}1_{I}(x))$
be a cut-off function for $I\subset I^*$ given by
Case 1. 
Let 
$\vp_0(x)=\Psi (R)^{-1}G^D_{\lambda}1_{B''}(x)$ where $B''=B(x_0,\delta R/2)$ 
and 
$\vp_1(x)=\vp_0(x)-\min_{y\in I^*}\vp (y)$, then by Lemma 3.3,
$$ \vp_1(x) \le c (s/R)^\th = L , \qq x\in I^*. $$
Let
$$\eqalign
{ A &=\int_{I} f^2 d\Gam(\vp,\vp), \cr
  D &= \int_{I^*} d\Gam(f,f)+ \Psi(s)^{-1} \int_{I^*}  f^2,\cr
  F &=\int_{I^*} f^2 \psi_s^2  d\Gam(\vp_1,\vp_1). \cr}$$

Now as
$$ d\Gam(f^2 \psi_s^2 \vp, \vp) \le  d\Gam(f^2 \psi_s^2 \vp_1 , \vp_0)
= f^2 \psi_s^2 d\Gam(\vp_1, \vp_0)  + \vp_1 d\Gam( f^2 \psi_s^2, \vp_0), $$
we have
$$ 
 A\le F= \int_{I^*} f^2 \psi_s^2 d\Gam(\vp_1,\vp_0) 
  =  \int_{I^*} d\Gam(f^2 \psi_s^2 \vp_1, \vp_0) 
  -  \int_{I^*} \vp_1 d\Gam(f^2 \psi_s^2 , \vp_0).   \eqno(3.8) $$
For the first term in (3.8)
$$  \eqalign{ 
 \int_{I^*} d\Gam(f^2 \psi_s^2 \vp_1, \vp_0)
  &=   \int_{X} d\Gam(f^2 \psi_s^2 \vp_1, \vp_0) \cr
 &= \sE_\lam(f^2 \psi_s^2 \vp_1,  \Psi (R)^{-1}G^D_{\lambda}1_{B''}) 
- \lam\int_X f^2 \psi_s^2\vp_1 \vp_0 d\mu\cr
 &\le  \sE_\lam(f^2 \psi_s^2 \vp_1, \Psi (R)^{-1}G^D_{\lambda}1_{B''}) 
=  \Psi (R)^{-1}
\int_{B''} f^2 \psi_s^2 \vp_1d\mu =0. \cr }$$
Here we used the fact that  $\vp_1\ge 0$ on $I^*$ and 
that the support of $\psi_s$ is in $I^*$, hence outside $B''$ (due to (3.7)).

The final term in (3.8) is handled, using the Leibniz and chain rules and 
Cauchy-Schwarz, as
$$ \eqalign{
  &\Bigl| \int_{I^*} \vp_1 d\Gam(f^2 \psi_s^2, \vp_0) \Bigr|  
\le 2\Bigl| \int_{I^*} \vp_1f\psi_s^2 d\Gam(f , \vp_0) \Bigr|  +
2\Bigl| \int_{I^*} \vp_1f^2\psi_s d\Gam( \psi_s, \vp_0) \Bigr|  
\cr
 &\le c \Big\{ \Big(\int_{I^*} \psi_s^2  d\Gam(f,f) \Big)^{1/2}
 + \Big(\int_{I^*} f^2  d\Gam(\psi_s,\psi_s){\Big)}^{1/2} \Big\} 
   \Big( \int_{I^*} \vp_1^2 f^2 \psi_s^2 d\Gam(\vp_0,\vp_0)\Big)^{1/2} \cr
 &\le  c D^{1/2} L  F^{1/2}, }$$ 
where we used Case 1 in the final line. 
Thus we obtain
$A\le F \le c DL^2$ so that (3.4) holds. 

\nl {\bf Case 3.} We finally consider the general case. 
When either (3.6) or (3.7) holds, the result is already proved in Case 2.
So assume that neither of them hold.
Then $I^*$ must intersect both $B(x_0,\delta R/2)$ 
and $B(x_0,\delta R)^c$, so $s\ge \delta R/4$.
We use Lemma 2.3 to cover 
$I$ with balls $B_i=B(x_i,c_1R)$, where $c_1\in (0, \delta/4)$
has been chosen small enough so that each 
$B_i^*:=B(x_i, c_1R/\delta)$  
satisfies at least one of (3.6) or (3.7). 
We can then apply (3.4) with $I$ replaced by each ball $B_i$:
writing $s'=c_1R$ we have
$$\int _{B_i} f^2 d\Gam(\vp,\vp) 
 \leq c_2(s'/R)^{2\th}  \Big(\int_{B^*_i} d\Gam(f,f)
+\Psi(s')^{-1} \int_{B^*_i} f^2d\mu\Big).$$
We then sum over $i$. Since no point
of $I^*$ is in more than $L_0$ (not depending on $x_0$ or $R$) 
of the $B^*_i$, and $s/c_1 \le s'\le s$, we obtain (3.4) for $I$. \qed

\ms \subsec{4. Sobolev inequalities and elliptic Harnack inequality}

\ms In this section we will prove Theorem 2.16 (a) $\Rightarrow$ (b). 
Assume that $X$ satisfies $\VD$, $\PI(\Psi)$ and $\CS(\Psi)$. 
Using Theorem 2.15 (c) $\Rightarrow$ (a)
and Lemma 2.7, it is enough to show $\VD+\PI(\Psi)+\CS(\Psi)
\Rightarrow \EHI$. 
For $x \in X$, $R\ge 0$ let $\vp=\vp_{x,R}$ be a cut-off function
given by $\CS(\Psi)$. We define the measure $\gam= \gam_{x,R}$ by
$$ d\gam = d\mu + \Psi(R) d\Gam(\vp,\vp). $$
We remark that we do not know if the measure $\gam$ satisfies volume doubling.
The first step in the argument is to use  $\CS(\Psi)$ to obtain a 
weighted Sobolev inequality. 
For any set $J\subset X$ set
$$J^s= \{ y: d(y,J)\le s\}.  $$

\proclaim Proposition 4.1.
Let $s\le R$ and $J\subset B(x_0,R)$ be a finite union of
balls of  radius $s$. 
There exist $\kappa>1$ and $c_1>0$ such that
$$ \Big(\mu(J)^{-1} \int_{J} |f|^{2\kappa} d\gamma\Big)^{1/\kappa} 
\le c_1 \Big(\Psi (R)\mu(J)^{-1} \int_{J^s} d\Gam(f,f) +   
  (s/R)^{-2\th}\mu(J)^{-1} \int_{J} f^2 d \gamma \Big). \eqno(4.1) $$

We omit the proof, since it is the same as that of Theorem 5.4 of [BB5].

\ms The next result is the generalization of  Lemma 4 of [M1] 
to the case of a MMD space.

\proclaim Lemma 4.2. Let $D$ be a domain in $X$,
let  $u$ be positive and harmonic in $D$, 
$v=u^k$, where $k \in \R$, $k\ne \fract12$,
and let $\eta$ be supported in $D$. 
Suppose $\int_D d\Gam(\eta,\eta)<\infty$, then
$$ \int_D \eta^2 d\Gam(v,v) \le
\Bigl(\frac{2k}{2k-1}\Bigr)^2 \int_D v^2 d\Gam(\eta,\eta) .$$

\proof Let $g \in \sF$ be supported by $D$. Then if $u'=G h$ where
$h=0$ on $D$ we have
$$ \int_D d\Gam(gu',u') = \int_X d\Gam(gu',u') =\int_X gu'h \,d\mu=0 .$$
Hence, approximating $u$ by functions of the form $u'$ we deduce that 
$$    \int_D d\Gam(gu,u)=0. $$ 
Using this, and taking $g=\eta^2 k^2 u^{2k-2}$, we conclude that
$$  \int_D \eta^2\, d\Gam(v,v) =   \int_D g\, d\Gam(u,u)
 = -  \int_D u \, d\Gam(g,u). \eqno(4.2)$$
Using the Leibniz and chain rules, the right hand side is equal to
$$ - 2k \int_D \eta v \, d\Gam(\eta,v)- (2k-2)\int_D \eta^2 d\Gam(v,v).$$
Thus, 
$$ \eqalign{\int_D \eta^2 d\Gam(v,v) 
&= -\frac {2k}{2k-1}\int_D\eta v\, d\Gam(v,\eta)\cr
&\le \frac{2|k|}{|2k-1|} \Big(\int_D \eta^2 d\Gam(v,v)\Big)^{1/2} 
\Big(\int_D v^2 d\Gam(\eta,\eta)\Big)^{1/2},\cr}$$ 
where we used Cauchy-Schwarz. Dividing and squaring, 
we obtain the result.\qed

\ms
Let $u$ be harmonic and nonnegative in $B(x_0,4 R)$.  By looking at
$u+\eps$ and letting $\eps\downarrow 0$ we may without loss of
generality suppose $u$ is strictly positive.  
Note that, as for a general MMD space we do not initially have any
{\sl a priori} continuity for $u$, we do not obtain a
pointwise bound in (4.3).

\proclaim Proposition 4.3. Let $v$ be either $u$ or $u^{-1}$.
There exists $c_1$ such that if $B(x,2r) \subset B(x_0,4R)$ and 
$0<q<2$, then
$${\es}_{B(x,  r/2)} v^{2q} \leq c_1 V(x,2r)^{-1}\int_{B(x,2r)}
\Big( \Psi(r) d\Gam(v^q,v^q) + v^{2q}\dm\Big) .\eqno(4.3)$$

\proof
Let $\vp_0$ be a (regularized) cut-off function given by $\CS(\Psi)$
for $B(x,r)$.
Let $h_n= 1-2^{-n}$, $0\le n\le \infty$, so that $0=h_0<h_\infty=1$.
For $k\geq 0$ set
$$ \vp_k(x)= (\vp_0(x)-h_k)^+, ~~~~
 d\gamma_0 =d\mu+ \Psi (r)d\Gam(\vp_0,\vp_0).$$
Set $A_k=\{x: \vp_0(x)>h_k\}$, and note that 
$B(x, r/2) \subset A_{n_0} \subset A_0 \subset B(x,r)$ for every $n_0$.
We therefore have, writing $V$ for $V(x,r)$,
$$ c_2 V \le \mu(A_k) \le V, \q k\ge 0. $$
The H\"older condition on $\vp_0$ given by $\CS(\Psi)$ 
implies that if $x\in A_{k+1}$ and $y\in A_k^c$, then
$d(x,y) \ge c_3 r 2^{-k/\th}$.  
Set $s_k = \half c_3 r 2^{-k/\th}$, and note that 
$\vp_k> c_4 2^{-k}$ on $A_{k+1}^{s_k}$.
Let $\{B_i\}$ be a cover of $A_{k+1}$ by balls of radius $s_k/2$, 
and let $J_{k+1} = \cup_i B_i$. Write $J'_{k+1} = J^{s_k/2}_{k+1}, 
A'_{k+1}=A_{k+1}^{s_k}$ and note that 
$A_{k+1} \subset J_{k+1} \subset J'_{k+1} \subset A'_{k+1}$.

{}From Proposition 4.1 with $f=v^p$ and $s$ replaced by $s_{k}/2$,
$$ \eqalignno{
 \Big( V^{-1} \int_{A_{k+1}} f^{2\kappa}d\gamma_0  \Big)^{1/\kappa} 
 &\le \Big( V^{-1} \int_{J_{k+1}} f^{2\kappa}d\gamma_0 \Big)^{1/\kappa} \cr
 &\le c_{5} V^{-1} \Big[ \Psi(r)  \int_{J'_{k+1}} d\Gam (f,f) + 
   (r/s_k)^{2\th} \int_{J'_{k+1}} f^2 d\gamma_0 \Big] \cr
 &\le c_6  V^{-1} \Big[ \Psi(r)  \int_{A'_{k+1}} d\Gam (f,f) + 
   2^{2k} \int_{A_k} f^2 d\gamma_0  \Big]. &(4.4)\cr }$$
By Lemma 4.2, we have the `converse to the Poincar\'e inequality'
for $f=v^p$, which controls the first term in (4.4).
$$\eqalign{ 
 \Psi(r)  \int_{A'_{k+1}} d\Gam (f,f)
 &\le \Psi(r) (c_7 2^{-k} )^{-2} \int_{A'_{k+1}} \vp_k^2 d\Gam (f,f) 
\le c_8 2^{2k} \Psi(r) \int_{A_{k}} \vp_k^2 d\Gam (f,f) \cr
&\le c_9 2^{2k} \Psi(r)
 \Big(\frac{2p}{2p-1}\Big)^2  \int_{A_{k}} f^2 d\Gam(\vp_k,\vp_k)
\le c_{10} 2^{2k} 
 \Big(\frac{2p}{2p-1}\Big)^2  \int_{A_{k}} f^2 d\gamma_0 . \cr }$$
We therefore deduce that

$$ \Big( V^{-1} \int_{A_{k+1}} f^{2\kappa}d\gamma_0 \Big)^{1/\kappa} 
\le c_{11} \Big(\frac{2p}{2p-1}\Big)^2 2^{2k}
     V^{-1} \int_{A_{k}} f^2 d\gamma_0. \eqno (4.5)$$

Now, similarly to the first part of Moser's argument [M1]
with $p_n=q \kappa^n$ for appropriate $q$, we have
$$\es_{B(x, r/2) }v    \le c_{12} \Bigl( V^{-1} \int_{B(x,r)} 
       v^{2q} d\gamma_0 \Bigr)^{1/(2q)}.$$
Using $\CS(\Psi)$ and $\VD$, we obtain (4.3) -- see [BB5]  Proposition 5.8 
for details. \qed
\ms
Recall that $\vp$ is a cut-off function for $B(x_0,R)$ given by $\CS(\Psi)$. 
We define 
$$ Q(t) = \{ x: \vp(x) > t \}, \q 0<t< 1, $$
and write $Q(1)$ for the interior of $\{x : \vp(x) \ge 1\}$.

\smallskip
Using Proposition 4.3, we have the following. 
See [BB5] Corollary 5.9 for the proof.  

\proclaim Corollary 4.4. Let $1 > s>t > 0$.
There exists $\zeta>2$ such that if $0<q< \frac13$,
$$ \es_{Q(s)} v^{2q} \le c_1 (s-t)^{-\zeta} 
   V(x_0,R)^{-1} \int_{Q(t)} v^{2q} d\gamma .\eqno (4.6)$$

Now our goal is to deduce the elliptic Harnack inequality. 
The following corresponds to the second part of Moser's arguments.

Let $w=\log u$, and write 
$\ol w= V(x_0,R)^{-1} \int_{B(x_0, R)} w\, d\mu$.

\proclaim Proposition 4.5. (cf. [BB5] Proposition 5.7, Corollary 5.10.)

\ni (a) There exists $c_1 $  such that 
$$ \int_{ B(x_0,2 R)} d\Gam (w,w)\le c_1 \frac{ V(x_0,R)}{\Psi (R)} .$$
(b) Let $1 \ge s>t >0$.
Then
$$ \int_{ \{|w-\ol w|> A \} \cap  Q(s)} d\gamma
  \leq c_2\frac{V(x_0,R) } {A^2}.$$

\proof
Again, this is essentially Moser's proof.
Let $\vp_1(x)$ be a cut-off function given by $\CS(\Psi)$
for the ball $B^*:=B(x_0,4R)$.  
So
$$ \int_{B(x_0,2 R)} d\Gam(w,w) \dmp \leq c
\int_{B^*} \vp_1^2 d\Gam(w,w) \dmp .$$
Applying (4.2) with $\eta=\vp_1$, $v=w$, $g=\vp_1^2/u^2$ and $D=B^*$,
we have 
$$\int_{B^*}\vp_1^2d\Gam (w,w)=-\int_{B^*}ud\Gam (\vp_1^2/u^2,u).$$ 
Using the Leibniz and chain rules, the right hand side is equal to
$$ - 2 \int_{B^*}\vp_1 d\Gam(\vp_1,w)+2\int_{B^*} \vp_1^2 d\Gam(w,w).$$
Thus, 
$$ \int_{B^*} \vp_1^2 d\Gam(w,w) = 2\int_{B^*}\vp_1 d\Gam(\vp_1,w)
\le 2 \Big(\int_{B^*} d\Gam(\vp_1,\vp_1)\Big)^{1/2} 
\Big(\int_{B^*} \vp_1^2 d\Gam(w,w)\Big)^{1/2},$$ 
where we used Cauchy-Schwarz. Dividing and squaring,  
$$ \int_{B^*} \vp_1^2d\Gam (w,w)
\leq 4\int_{B^*} d\Gam (\vp_1,\vp_1). $$
Finally, using $\CS(\Psi)$ in $B^*$ with 
$f\in \sF$ such that $f|_{B(x_0,8R)}\equiv 1$ (since $(\sE,\sF)$ 
is a regular Dirichlet form, such an $f$ exists) 
and $\VD$  we deduce that 
$$ \int_{B^*} d\Gam (\vp_1,\vp_1)
 \le c  \Psi (R)^{-1} V(x_0,R). $$
The proof of (b) is the same as that of [BB5] Corollary 5.10, 
so we omit it.\qed

In order to get the Harnack inequality the argument in [M2] required a 
generalization of the John-Nirenberg inequality with a complicated proof.
Bombieri [Bom] found a way to avoid such an argument
for elliptic second order differential equations. Moser (Lemma 3 in [M3]) 
carried the idea over to the parabolic case and Bombieri and Giusti (Theorem 4 
in [BG]) obtained the inequality in an abstract setting. (See also Lemma 2.2.6 in [SC2].) 
This argument can be applied to our setting (with suitable 
modifications) and we can show that Corollary 4.4 and Proposition 4.5 (b) give
$$\es_{B(x_0,R/2) } \log u \leq c_1,$$
(see [BB5] Lemma 5.11 and Theorem 5.12 for a detailed proof).
Let $v=u^{-1}$. The same argument implies 
$\es_{ B(x_0,R/2)}  \log v \leq {c_1}$,
or $ \ei_{B(x_0,R/2) }  \log u \geq -{c_1}.$
Combining we deduce
$$ e^{-c_1} \le \ei_{ B(x_0,R/2) } u 
\le \es_{B(x_0,R/2) } u \le e^{c_1}.$$
We thus obtain the following. 
\proclaim Theorem 4.6. There exists $c_1$ such that if $u$ is nonnegative
and harmonic in $B(x_0,4R)$, 
then $$ {\es_{B(x_0, R/2)} u} \le c_1 { \ei_{B(x_0, R/2)} u} .$$ 

\longproof{of Theorem 2.16 (a) $\Rightarrow$ (b)}

\ni
As we mentioned in the beginning of this section, using Theorem 2.15 (c) 
$\Rightarrow$ (a) and Lemma 2.7, it is enough to show $\VD+\PI(\Psi)+\CS(\Psi)
\Rightarrow \EHI$. 

Let $u$ be nonnegative and harmonic in $B(x_0,4R)$. 
Suppose $x_1$ and $r$ are such that $B(x_1,3r)\subset B(x_0,4R)$.
By looking at $Cu+D$ for suitable constants $C$ and $D$, we may suppose that
$\es_{B(x_1,2r)} u=1$ and $\ei_{B(x_1,2r)} u=0$. 
Hence by Theorem 4.6 
we have
$$ \es_{B(x_1,r)} u-\ei_{B(x_1,r)} u \le (1-c_1^{-1}) \es_{B(x_1,r)} u
 \le   (1-c_1^{-1}). $$
So if $\rho=1-c_1^{-1}$ then 
$$\es_{B(x_1,r)} u-\ei_{B(x_1,r)} u \leq \rho
\big[\es_{B(x_1,2r)} u-\ei_{B(x_1,2r)}\big].$$
It follows easily that 
$$\es_{B(x_1,r)} u-\ei_{B(x_1,r)} u \leq c_2 r^\gamma \eqno (4.7)$$
for some $\gamma>0$. Define $\wh u(x_1)=\lim_{r\to 0} \es_{B(x_1,r)} u$.
If one takes a countable basis $\{B_i\}$ for $X$ and excludes those
points  $x\in B_i$ such that $u(x) \notin [\ei_{B_i} u, \es_{B_i} u]$,
then for every other $x$ it is easy to see, using (4.7), that $u(x)=\wh u(x)$. 
Thus, $\wh u$ is equal to $u$ for $\mu$-almost every $x$.
Moreover, from (4.7) we see that $\wh u$ is H\"older continuous. Recall that
in our definition of harmonic function we take a quasi-continuous
modification as defined in [FOT]. We conclude $u=\wh u$ quasi-everywhere,
and so $u$ has a quasi-continuous modification that is continuous.
Using this modification and  Theorem 4.6, we have 
$$ {\sup_{B(x_0, R/2)} u} \le c_1 { \inf_{B(x_0, R/2)} u} .$$ 
The elliptic Harnack inequality $\EHI$  now follows 
by a covering argument. \qed

\bigb {\bf 5. Stability under rough isometries. }

\ms In this section we prove Theorem 2.21.
It is known (and easy) that the condition $\VD$ is  
stable under rough isometry, see [CS, HK2] etc.

\ms In [HK2] the stability of $\VD+ \PI(\Psi)$ and $\VD+ \CS(\Psi)$ 
under rough isometries is proved in the case when the spaces are 
graphs with controlled weights. Since  rough isometry is an equivalence 
relation, to prove Theorem 2.21 it
is enough to prove that if $X$ is a MMD space satisfying $\VD_\loc$,  
$\PI(\betaloc)_\loc$ and $\CS(\betaloc)_\loc$  
and $G$ is a graph constructed by taking an appropriate net
of $X$, (so that $X$ and $G$ are  roughly isometric), then
$\VD$, $\PI(\Psi)$ and $\CS(\Psi)$ hold for $X$ if and only if 
they hold for $G$.

\ms Let $X$ be a MM space satisfying $\VD_\loc$. Let $G\subset X$
be a maximal set such that
$$ d(x,y)\ge 1 \hbox{ for } x,y \in G, x\neq y. $$
Thus $B(x,\half)$, $x\in G$ are disjoint, and $\cup_{x\in G} B(x,1) = X$. 
Give $G$ a graph structure by letting $x\sim y$ if $d(x,y)\le 3$.
Let $\dg$ be the usual graph distance on $G$, and write
$B_G(x,r)$ for balls in $G$.
It is straightforward to check that $G$ is connected, and that
$$ \fract13 d(x,y) \le d_G(x,y) \le d(x,y) +1, \q x,y \in G.$$ 
Since $X$ satisfies $\VD_\loc$ we have, as in Lemma 2.3 of [Kan1], that
the vertex degree in $G$ is uniformly bounded. 

For each $x\sim y$ in $G$ let $z_{xy}$ be the midpoint of
a geodesic connecting $x$ and $y$, and $A_{xy}=B(z_{xy},5/2)$,
so that $B(x,1) \subset A_{xy} \subset B(x,4)$.
Let $ \nu_{xy}=0$ if $x \not\sim y$, and if  $x \sim y$ let
$$ \nu_{xy} = \mu(A_{xy}). $$
As usual we set $\nu_x = \sum_{y\sim x} \nu_{xy}$. 
Write $A_x =\cup_{y\sim x} A_{xy}$.
Since $X$ satisfies $\VD_\loc$, we have
$$ \mu(B(x,1)) \le \nu_x  \le c_1 \mu(A_x)
 \le c_1 \mu(B(x,4)) \le c_2  \mu(B(x,1)), \eqno(5.1) $$
and using (2.3) it is easy to verify that $(G,\nu)$ has controlled
weights.

Define  $\imath: G \to X$ by  $\imath(x)=x$. We have

\proclaim Proposition 5.1. Let $X$ be a MM space satisfying
$\VD_\loc$. Then the associated weighted graph $(G,\nu)$ has
controlled weights and $\imath$ is a rough isometry.

\ni In the following, we abuse notation and denote 
the image of $\imath$ by the same character as its pre-image.

\ms To prove Theorem 2.21(b) we will need to transfer functions between 
$C(G, \bR_+)$ and $C(X,\bR_+)$. 
Let $f\in C(X,\bR_+)$. Define
$$ \wh f (x) =  \mu(B(x,1))^{-1} \int_{B(x,1)} f d\mu, \q x\in G. \eqno(5.2)$$
 The transfer in the other direction requires a bit more care. 
Using $\VD_\loc$ and $\CC_\loc$, 
we can find a partition
of unity $\{\psi_x\}_{x\in G}$ where each $\psi_x: X\to \bR$ 
is quasi continuous and satisfies the following:
\ms
\i{(i)} $\psi_x(w) =1$ for $w\in B(x,\fract14)$,
\i{(ii)} $\psi_x(w) =0$ for $w\in B(x,\fract32)^c$,
\i{(iii)} $ \int_{B(x,1)} d\Gam(\psi_x,\psi_x) \le c V(x,1).$

\ni If we assume $\CS(\bar\beta)_\loc$, then we can choose 
$\{\psi_x\}_{x\in G}$ that further satisfies the following:
\ms
\i{(iii')} For each $z\in G$, $s\le 1$ and $f\in \sF$, 
$$ \int_{B(z,s)} f^2 d\Gam(\psi_x,\psi_x) \le 
 c_1 s^{2\th} \Big( \int_{B(z,2s)} d\Gam(f,f) 
+ \Psi(s)^{-1} \int_{B(z,2s)} f^2  d\mu\Big).\eqno (5.3)$$

\ms
Now if $g: G \to \bR_+$ set 
$$ \wt g(z) =\sum_{x \in G} g(x) \psi_x(z). \eqno(5.4)$$

Note that $\wt g: X\to \bR$ is quasi continuous.
Set also, if $f: G \to \bR$, $k \in \bN$,
$$ V_k f(x) = \sup_{z:d_G(z,x)\le k} |f(x)- f(z)|. $$

\proclaim Lemma 5.2. Let $X$ be a MMD space satisfying
$\VD_\loc$ and $\CC_\loc$. 
Let  $f: G \to \bR_+$, and $x\in G$.
\nl (a) If $\psi_z(w)>0$ for some $w\in B(x,1)$ 
then $d(x,z)<3$ and $x\sim z$.
\nl (b) If $\psi_z(w)>0$ for some $w\in A_x$ then $d_G(x,z) \le 4$, and
$$ |f(x)-\wt f(w)| \le V_4 f(x), \q w\in A_x. $$
(c) Let $A \subset G$, and $A'=\{y: d_G(y,A)\le 4\}$. Then
$$ \sum_{z \in A} V_4 f(z)^2 \nu_z 
  \le c_1 \sum_{y,z \in A'} (f(y)-f(z))^2 \nu_{yz}. $$
(d) 
$$  c_2\sum_{y \in G\cap B(x,r-1)} f(y)^2 \nu_y
\le \int_{B(x,r)} \wt f(w)^2 d\mu(w) 
 \le c_3\sum_{y \in G\cap B(x,r+2)} f(y)^2 \nu_y. $$
(e) $$\int_{B(x,1)} d\Gam(\wt f,\wt f) \le c_4V_1f(x)^2V(x,1).$$
(f) If $X$ further satisfies $\CS(\bar\beta)_\loc$. Then, 
for each $s\le 1$ and $h\in\sF$,  
$$ \int_{B(x,s)} h^2 d\Gam(\wt f,\wt f) \le 
 c_5 s^{2\th} V_1 f(x)^2 \Big( \int_{B(x,2s)} d\Gam(h,h) 
+ \Psi(s)^{-1} \int_{B(x,2s)} h^2  d\mu\Big).$$

\proof The proof of (a)-(c) and the second inequality of 
(d) is simple, and is the same as that of [B2] Lemma 5.6 (a)-(d).
So we will give the rest of the proof.

For the first inequality of (d), there is nothing to prove when 
$r<1$, so assume $r\ge 1$. Using (5.1) and $\VD_\loc$, we see that 
that the left hand side of (d) is bounded by
$$
c\sum_{y \in G\cap B(x,r-1)}f(y)^2V(y,1/4)=c\sum_{y \in G\cap B(x,r-1)} 
\int_{B(y,1/4)}\wt f(w)^2d\mu (w),$$
where the equality holds because for each $y\in G$ and 
$w\in B(y,1/4)$, $\wt f(w)=f(y)$.
Since $\{B(y,1/4)\}_{y\in G}$ are disjoint, the right hand side is bounded
from above by $\int_{B(x,r)} \wt f(w)^2 d\mu(w)$.

Before proving (e), let us prove (f). 
By (a) we can write, for $w\in B(x,1)$,
$$ \wt f(w) = \wt f(x) + \sum_{z\sim x} \psi_z(w) (f(z)- f(x)). $$
Hence
$$ \eqalignno{
  \int_{B(x,s)} h^2d\Gam( \wt f, \wt f) &= 
\sum_{z\sim x} \sum_{z' \sim x} 
 (f(z)- f(x))(f(z')- f(x))  \int_{B(x,s)} h^2d\Gam( \psi_z, \psi_{z'}) \cr
 &\le c \sum_{z\sim x}  (f(z)- f(x))^2   \int_{B(x,s)} h^2
d\Gam( \psi_z, \psi_{z}) \cr
 &\le c \sum_{z\sim x}  V_1f(x)^2  \int_{B(x,s)} h^2 d\Gam( \psi_z, \psi_{z}) 
&(5.5)\cr
&\le c' s^{2\th} V_1 f(x)^2 \Big( \int_{B(x,2s)} d\Gam(h,h) 
+ \Psi(s)^{-1} \int_{B(x,2s)} h^2  d\mu\Big), \cr}$$
where we used (5.3) in the last inequality.

Now (e) is easy. Using (5.5) with $s=1$ and $h|_{B(x,2s)}\equiv 1$, we have 
$$ \int_{B(x,1)}d\Gam( \wt f, \wt f)\le 
c \sum_{z\sim x}  V_1f(x)^2  \int_{B(x,1)}d\Gam( \psi_z, \psi_{z}) 
\le c' V_1f(x)^2 V(x,1),$$
where $\CC_\loc$ and $VD_\loc$ are used for the second inequality.
\qed

\proclaim Lemma 5.3.  Let $X$ be a MMD space satisfying
$\VD_\loc$ and $\PI_\loc$. 
Let $g: X \to \bR_+$,  $x\in G$, and $y\sim x$. Then
$$  (\wh g(x)-\wh g(y))^2 \nu_{xy} \le c \int_{A_{xy}} d\Gam(g,g). $$

The proof is simple, and it is the same as that of [B2] Lemma 5.7, so we omit it. 
\ms
In the arguments that follow, we will use the fact, given in
Remarks 2.6, that to verify $\CS(\Psi)$ it is enough to
do so for any $\delta>0$ in (2.6) and $\lam>0$ in (2.7). 

\proclaim Proposition 5.4.  Let $X$ be a MMD space satisfying
$\VD_\loc$, $\PI_\loc$ and $\CC_\loc$, and $(G,\nu)$ be the 
associated weighted graph.
\nl (a) Suppose that $X$ satisfies $\VD$ and 
$\PI(\Psi)$. Then $(G,\nu)$  satisfies $\VD$ and $\PI(\Psi)$.
\nl (b) Suppose that $X$ satisfies $\VD$ and $\CS(\Psi)$. 
Then $(G,\nu)$ satisfies $\VD$ and $\CS(\Psi)$.

\proof
As mentioned above, $\VD$ is preserved, so it is enough to prove $\PI(\Psi)$ 
and $\CS(\Psi)$ respectively. 

We first prove (a). 
Let $x\in G$, $R>0$, $f: G\to \bR$ and let
$B_G(x,R)$ be a ball in $G$. It is enough to 
prove the following weak Poincar\'e inequality;
$$\sum_{y\in B_G(x,R)}(f(y)-\bar f_{B_G(x,R)})^2\nu_x\le 
c\Psi(R)\sum_{y\in B_G(x,c'R)}\Gamma_G(f,f)(y),\eqno(5.6)$$
where $c>0, c'>1$ are constants.
Indeed, an argument such as that in Jerison [Je] (see also [HaKo] for a more general 
formulation) can be used to derive the (strong) Poincar\'e inequality $\PI(\Psi)$ 
from $\VD$ and (5.6).

We note that for fixed $c_0\ge 1$, $\PI(\Psi)$ for $R\le c_0$ 
always holds on $(G, \nu)$.
Indeed, if we let $M_{x}:=\max\{|f(x')-f(y')|: x',y'\in B_G(x_0, R)\}$, 
then, since $\min_{x\in B} f(x)\le {\bar f}_B \le \max_{x\in B} f(x)$, we have
that the left hand side of (5.6) is bounded by
$$  2M_{x}\sum_{y\in B_G(x, R)}\nu_y
\le C_{c_0}\sum_{y\in B_G(x, R)}\Gam_G(f,f)(y),$$
where $C_{c_0}>0$ depends on $c_0$. 
This is turn is less than or equal to the right hand side of (5.6).
So it is enough to prove (5.6) for $R\ge c_0$.
Applying Lemma 5.2.(d) for $f-\alpha$ where $\alpha=V(x,R)^{-1}
\int_{B(x,R)}\wt fd\mu$ and noting $B_G(x,R')\subset G\cap B(x,R-1)$
where $R'=(R-1)/3$, we have 
$$\int_{B(x,R)} (\wt f(w)-\alpha)^2 d\mu(w)\ge 
c\sum_{y \in B_G(x,R')} (f(y)-\alpha)^2 \nu_y
\ge c\sum_{y \in B_G(x,R')} (f(y)-\bar f_{B_G(x,R')})^2 \nu_y.$$
 Here the last inequality is because the minimum of the middle term 
(as a function of $\alpha$) is attained when $\alpha= \bar f_{B_G(x,R')}$. 
On the other hand, Using Lemma 5.2.(e), we have 
$$ \eqalign{
\int_{B(x,R)}d\Gam (\wt f,\wt f)&\le \sum_{y\in B_G(x,R)}\int_{B(y,1)}
d\Gam (\wt f,\wt f)\le c\sum_yV_1f(y)^2V(y,1)\cr
&\le c'\sum_{y\in B_G(x,R+2)}\Gam_G(f,f)(y).}$$
Combining the estimates with $\PI(\Psi)$ for $X$ gives (5.6) for $R\ge c_0$.

We next prove (b). We need to construct a
cut-off function $\wh \vp$ satisfying (a)--(d) of Definition 2.5.
If $R\le c_0$ then  it is easy to check that we can take 
$\wh\vp(x)$ to be the indicator of $B_G(x_0,R/2)$. 

So assume $R>c_0$. We can find a constant $c_1$ such that
$$ B_G(x_0, c_1 R) \subset G\cap B(x_0,R/8-6) \subset  G\cap B(x_0,R/4+6)
\subset B_G(x_0,R). $$
It is enough to construct a cut-off function $\wh\vp$ for 
$ B_G(x_0, c_1 R) \subset  B_G(x_0,R)$.
Let $\vp$ be a cut-off function for $B(x_0,R/8) \subset B(x_0,R/4)$,
and let $\wvp$ be given by (5.2). Properties (a)--(c) of Definition
2.5 are easily checked, and it remains to verify the weighted
Sobolev inequality (2.5). 

Let $x_1\in G$, $1\le s\le R$, and $A_G=B_G(x_1,s)$.
Choose $c_2, c_3$ so that
$$ A_G \subset B(x_1,c_2 s-6)\cap G \subset B(x_1,2 c_2 s)\cap G 
 \subset B_G(x_1,c_3s-6 ). $$
Write $A'_G=B_G(x_1, c_3 s)$, and let  $f: A'_G \to \bR_+$. 
We extend $f$ to $G$ by taking $f$ to be zero outside $A'_G$, and
define $\wt f$  by (5.4).

Let $x\in G$, and $y\sim x$. Then by Lemma 5.2(b) and Lemma 5.3
$$ \eqalign{
 f(x)^2  (\wvp(x)-\wvp(y))^2 \nu_{xy} 
 &\le c  \int_{A_{xy}}  f(x)^2 d\Gam(\vp,\vp)(w) \cr
 &\le 2 c  \int_{A_{xy}}  \wt f(w)^2 d\Gam(\vp,\vp)(w) 
  +  2 c  \int_{A_{xy}} V_4 f(x)^2 d\Gam(\vp,\vp)(w). \cr} $$
Therefore 
$$ \eqalignno{
 \sum_{x \in A_G}  \sum_{y\sim x} &f(x)^2 (\wvp(x)-\wvp(y))^2 \nu_{xy} \cr
 &\le c \sum_{x \in A_G} \sum_{y\sim x}   
  \int_{A_{xy}}  \wt f(w)^2 d\Gam(\vp,\vp)(w) 
   + c  \sum_{x \in A_G} \sum_{y\sim x} 
   \int_{A_{xy}} V_4 f(x)^2 d\Gam(\vp,\vp)(w) \cr
 &\le c   \int_{B(x_1,c_2 s)}  \wt f(w)^2 d\Gam(\vp,\vp)(w)
 + c  \sum_{x \in A_G}  \sum_{y\sim x} 
    V_4 f(x)^2 \int_{A_{xy}} d\Gam(\vp,\vp)(w) . &(5.7) \cr}$$
Applying $\CS(\Psi)$ to $\vp$ in the ball $A_{xy}$ gives
$$  \int_{A_{xy}} d\Gam(\vp,\vp) \le c R^{-2\th} \mu(B(z_{xy},5))
 \le c'  R^{-2\th} \nu_{xy}. \eqno(5.8) $$
Therefore, using Lemma 5.2(c), the second term in (5.7) is bounded by
$$ 
  c R^{-2\th} \sum_{x \in A_G}  \sum_{y\sim x}  V_4 f(x)^2 \nu_{xy} \
 \le c R^{-2\th} \sum_{x\in A'_G} \Gam_G(f,f)(x). \eqno(5.9) $$

Using (2.5) for $\vp$ gives 
$$ \eqalignno{
 \int_{B(x_1,c_2s )}  &\wt f(w)^2 d\Gam(\vp,\vp)(w) \cr
& \le c (s/R)^{2\th} 
 \Big( \int_{B(x_1,2c_2s)}  d\Gam(\wt f,\wt f)
 + \Psi(s)^{-1}  \int_{B(x_1,2c_2s)} \wt f^2 d\mu \Big). &(5.10) \cr}$$
By Lemma 5.2(e), 
$$ \eqalignno{
  \int_{B(x_1,2c_2s)} d\Gam(\wt f,\wt f)
  &\le \sum_{x\in G\cap B(x_1,2c_2s+1 )}  \int_{B(x,1)} 
d\Gam(\wt f,\wt f)
\le c\sum_{x\in G\cap B(x_1,2c_2s+1 )} V_1 f(x)^2 V(x,1)\cr
 &\le c'\sum_{x,y \in A'_G} (f(x)-f(y))^2  \nu_{xy}, &(5.11) \cr} $$
while by Lemma 5.2(d) 
$$  \int_{B(x_1,2c_2s)} \wt f(w)^2 d\mu(w)
  \le  \sum_{x\in  B_G(x_1, c_3 s)} f(x)^2 \nu_x. \eqno(5.12)$$
Combining the estimates (5.7)--(5.12) completes the proof. \qed

\proclaim Proposition 5.5. Let $X$ be a MMD space satisfying
$\VD_\loc$, $\PI_\loc$,  and $(G,\nu)$ be the 
associated weighted graph.
\nl (a) Suppose that $X$ satisfies $\PI(\betaloc)_\loc$. If 
$(G, \nu)$  satisfies $VD$ and $\PI(\Psi)$ 
then $X$ satisfies $VD$ and $\PI(\Psi)$.  
\nl (b) Suppose that $X$ satisfies $\CS(\betaloc)_\loc$. If 
$(G, \nu)$  satisfies $VD$ and $\CS(\Psi)$ 
then $X$ satisfies $VD$ and $\CS(\Psi)$.  

\proof We first prove (a). 
Let $x_0\in X$, $R>0$ and 
$f\in \sF$. As mentioned in the proof of Proposition 5.4.(a),
it is enough to prove the following weak Poincar\'e inequality;
$$\int_{B(x_0,R)}(f(y)-\bar f_{B(x_0,R)})^2\mu (y)\le 
c\Psi(R)\int_{B(x_0,c'R)}d\Gamma(f,f),\eqno(5.13)$$
where $c>0, c'>1$ are constants. When $R\le 1$, this can be obtained 
from $\PI(\betaloc)_\loc$, so assume $R>1$. 
Using $\PI(\betaloc)_\loc$ with $R=1$, we have for each $x\in G$,
$$\int_{B(x,1)}f^2d\mu- \wh f(x)^2\nu_x=\int_{B(x,1)}(f-{\wh f}~)^2d\mu 
\le c\int_{B(x,1)}d\Gam (f,f).$$
Summing this over $x\in B(x_0,R)\cap G$, we have 
$$\int_{B(x_0,R)}f^2d\mu\le c\Big(\sum_{x\in B(x_0,R)\cap G}\wh f(x)^2\nu_x+
\int_{B(x_0,R+1)}d\Gam (f,f)\Big).$$
Putting $f-\alpha$ instead of $f$ where $\alpha=\sum_{x\in B(x_0,R)\cap G}\wh f(x)\nu_x
/(\sum_{x\in B(x_0,R)\cap G}\nu_x)$ and using the fact that 
$\int_{B(x_0,R)}(f-\beta)^2d\mu$
(as a function of $\beta$) attains its minimum when $\beta=\bar f_{B(x_0,R)}$, we obtain  
$$\int_{B(x_0,R)}(f-\bar f_{B(x_0,R)})^2d\mu\le c\Big(\sum_{x\in B(x_0,R+1)\cap G} 
(\wh f(x)-\alpha)^2\nu_x+\int_{B(x_0,R+1)}d\Gam (f,f)\Big).$$
On the other hand, using Lemma 5.3. and summing, we have 
$$\sum_{x\in B(x_0,R+1)\cap G}\Gam_G(\wh f~,\wh f~)(x)\le 
c \int_{B(x_0,R+5)}d\Gam (f,f).$$ 
Combining the estimates with $\PI(\Psi)$ for $G$ gives (5.13) for $R\ge 1$.

We next prove (b). 
Let $B=B(x_0,R)$ be a ball in $X$. 
If $R\le c_1$ then we can use $\CS(\bar\beta)_\loc$ to 
construct a cut-off function $\vp$ for $B$. So assume  $R\ge c_1$.
We can therefore assume that $x_0\in G$. 

Given $A\subset G$ write $A^{(1)}=\cup_{x\in A} B(x,1)$. 
We can find $c_i$ such that 
$$ B(x_0, c_1 R) \subset B_G(x_0, c_2R-6)^{(1)}
 \subset  B_G(x_0, 2c_2R+6)^{(1)} \subset  B(x_0,  R). $$
Let $\vp_G$ be a cut-off function for 
$B_G(x_0, c_2R) \subset  B_G(x_0, 2c_2R)$,
and let 
$$ \vp(w) = \wt \vp_G(w) =  \sum_{z\in G} \vp_G(x)\psi_z(w). $$
Properties (a)--(c) of $\vp$ follow easily from those of $\vp_G$, 
and it remains to verify (2.5). 

Let $B_1=B(x_1,s)$ with $s\in (0,R)$.
If $s\le c_3$ then, applying Lemma 5.2(f) and noting 
$V_1 \vp(x)\le c R^{-2\th}$, we have  
$$ \int_{B(x_1,s)} g^2 d\Gam(\vp,\vp) 
  \le c(s/R)^{2\th}\Big( \int_{B(x_1,2s)} d\Gam(g,g) 
+ \Psi(s)^{-1} \int_{B(x_1,2s)} g^2  d\mu\Big). $$

Now suppose $s\ge c_3$. Then we can assume $x_1\in G$,
and there exist $c_i$ so that
$$ B(x_1,s) \subset B_G(x_1, c_4s-6)^{(1)} \subset 
 B_G(x_1, 2 c_4s+6)^{(1)} \subset  B(x_1,c_5 s-6). $$
Let  $g:B(x_1,c_5 s)\to \bR_+$.  Define $\wh g$ on 
$ B_G(x_1, 2 c_4s+6)$ by (5.2). 
Then   
$$ \eqalignno{
 \int_{B(x_1,s)} g^2 d\Gam(\vp,\vp)
 &\le \sum_{x \in  B_G(x_1, c_4s)} \int_{B(x,1)} g(w)^2  d\Gam(\vp,\vp)(w)\cr
 &\le 2 \sum_{x \in  B_G(x_1, c_4s)} 
     \int_{B(x,1)} (g(w)-\wh g(x))^2 d\Gam(\vp,\vp)(w) &(5.14)\cr
 &\q + 2\sum_{x \in  B_G(x_1, c_4s)} 
    \int_{B(x,1)} \wh g(x)^2 d\Gam(\vp,\vp)(w).  &(5.15) \cr }$$
By Lemma 5.2(e) 
with $s=1$, the term (5.14) is bounded by
$$ c R^{-2\th}  \sum_{x \in  B_G(x_1, c_4s)} 
  \Big(\int_{B(x,1)}  (g(w)-\wh g(x))^2 d\mu
+\int_{B(x,1)}d\Gam (g,g)\Big), $$
and using $\PI_\loc$ this is bounded by
$$ c' R^{-2\th}  \sum_{x \in  B_G(x_1, c_4s)} \int_{B(x,1)} d\Gam(g,g)
  \le  c'' R^{-2\th}  \int_{B(x_1,c_5 s)} d\Gam(g,g). \eqno(5.16) $$
For the term (5.15), 
by Lemma 5.2(e) 
and (2.11) for $\vp_G$,
$$ \eqalign{ 
 \sum_{x \in  B_G(x_1, c_4s)}& \wh g(x)^2 \int_{B(x,1)} d\Gam(\vp,\vp)(w) \cr
& \le \sum_{x\in  B_G(x_1, c_4s)} 
  \wh g(x)^2 V_1 \vp_G(x)^2 V(x,1) \cr
&\le  c \sum_{x\in  B_G(x_1, c_4s)} 
  \wh g(x)^2 \Gam_G(\vp_G,\vp_G)(x) \cr
 &\le c(s/R)^{2\th} \Big( 
  \sum_{x\in  B_G(x_1, 2c_4s)} \Gam_G(\wh g,\wh g)(x)
 + \Psi(s)^{-1}  \sum_{x \in  B_G(x_1, 2c_4s)} \wh g(x)^2 \nu_x \Big).
}$$
Using Lemma 5.3 for the first term, and an easy bound for the second,
(2.5) now follows. \qed

We finally mention several applications of 
our results.

\ms\ni 1. {\sl  Dirichlet forms on generalized Sierpinski carpets}

\ni In [BB1, BB2] a `nice' diffusion 
process on a generalized Sierpinski carpet $F$ is constructed and it is
proved that its heat kernel satisfies $\HK(\Psi)$ with
$\Psi(s)=s^{d_w}$ for some $d_w\ge 2$.  Let the corresponding
Dirichlet form be denoted by $\sE_1$.  On the other hand, a
self-similar Dirichlet form  $\sE_2$ was constructed on the carpet by an
averaging method in [KZ] 
(see also [BB2] Remark 5.11 and [HKKZ]).
Because of the possible lack of
uniqueness for the `nice' diffusion, it is not known if the
corresponding diffusions coincide.  
Further, it was not clear that 
the heat kernel corresponding to $\sE_2$ satisfies $\HK(\Psi)$. 
(If one tries to apply the coupling methods used in [BB2], this
works for a dense set of starting points, but not for all starting
points because one does not know {\sl a priori} the continuity of
harmonic functions.) 
\proclaim Theorem 5.6. Let $(F,d)$ be a  generalized Sierpinski carpet
with walk dimension $d_w$, 
and  $\sE_2$ be the Dirichlet form constructed in  [KZ]. Then
the MMD space $(F,d, \sE_2)$ satisfies $\HK(d_w)$ and 
$\PHI(d_w)$. 

\proof 
It is easy to check that $\sE_2$ is a bounded perturbation of $\sE_1$
in the sense of (2.12) (see [Hin] Section 5.2 and [KZ] Section 6). 
Thus $\HK(\Psi)$ for $\sE_2$ holds by Theorem 2.19. \qed

In [Hin], general criteria are given for energy measures of
self-similar Dirichlet forms on self-similar sets to be singular with
respect to Bernoulli type measures. Roughly speaking, the main
theorems imply that the energy measures are singular with respect to
the Hausdorff measure if the elliptic Harnack inequality holds and the
walk dimension of the corresponding process (which corresponds to
$\beta$ and $\betaloc$ in this paper) is greater than $2$.  As a
consequence of the last paragraph, these conditions hold for the
self-similar Dirichlet forms on the Sierpinski carpets and we can
conclude that the energy measures of the Dirichlet forms (both for
$\sE_1$ and $\sE_2$ above) are singular with respect to Hausdorff
measure. (In fact, we can also verify (O) in [Hin], so using Theorem
2.4 in [Hin] we can further conclude the singularity with respect to a
certain class of Bernoulli type measures.) Note that this was not
proved in [Hin] for higher dimensional carpets 
since at that point the elliptic Harnack inequality had not been proved for 
the self-similar Dirichlet form $\sE_2$.

\ms\ni 2. {\sl Weighted graphs and manifolds associated with MM spaces 
satisfying $\VD_\loc$}

\ni 
In this section, we have explained how to construct weighted graphs from 
MM spaces satisfying $\VD_\loc$. There is a natural way to construct `jungle gym' type
manifolds from the weighted graphs (see [BCG, Kan2, PS] etc.). Roughly, this can be done 
by replacing the edges by tubes of length $1$, and by gluing the tubes
together smoothly at the vertices (see the right side of Figure 1).  
By Theorem 2.21, we know that if the original MM space satisfies $\PHI(\Psi)$ 
for some $\Psi(s)=s^{\betaloc}\vee s^{\beta}$, then the network on the associated
weighted graph and the Laplacian on the associated manifold satisfy $\PHI(\Psi')$ with 
$\Psi'(s)=s^{2}\vee s^{\beta}$. Further, any uniformly elliptic
operator in  divergence form on a 
manifold which is roughly isometric to the MMD space satisfies $\PHI(\Psi')$.

This fact is useful in fractal contexts, since $\PHI(\Psi)$ with $\Psi(s)=s^{d_w}$ 
for some $d_w\ge 2$ is proved for various `regular' fractals. Our result thus gives 
an alternative proof of the results in [BB3, BB4, Jo] and the heat kernel results in 
[HK2]. 

\medskip
\noindent {\bf Acknowledgment. }
The authors thank 
Zhenqing Chen, Alexander Grigor'yan, Masanori Hino, 
Laurent Saloff-Coste and Karl-Theodor Sturm 
for valuable comments, 
and Masayoshi Takeda for informing us about the references [LJ, Mos]. 

\def \nr { \nl }
\bigbreak \centerline {\bf References}
\bigskip
\nr [A]  D.G. Aronson. Bounds on the fundamental solution of a 
parabolic equation. {\it Bull.\ Amer.\ Math.\ Soc.\ \bf 73} (1967), 890--896.
\nr
[B1]
M.T. Barlow. Diffusions on fractals. In:
{\sl Lectures on Probability Theory and Statistics,
 Ecole d'\'Ete de Probabilit\'es de Saint-Flour XXV - 1995}, 1--121, 
Lect. Notes Math. {\bf 1690}, Springer, 1998.  
\nr
[B2] M.T. Barlow. 
Anomalous diffusion and stability of Harnack inequalities.  
To appear {\it Surveys in Differential Geometry IX}. 
\nr
[BB1] M.T. Barlow, R.F. Bass. 
Transition densities for Brownian motion on the Sierpinski carpet.
{\it Probab. Theory Rel. Fields} {\bf 91} (1992), 307--330.
\nr
[BB2] M.T. Barlow, R.F. Bass. Brownian motion and harmonic analysis on
Sierpinski carpets. {\it Canad. J. Math.} {\bf 51} (1999), 673--744.
\nr 
[BB3] M.T. Barlow, R.F. Bass. Random walks on graphical Sierpinski 
carpets. In: {\sl Random walks and discrete potential theory}  
(Cortona, 1997), 26--55, {\sl Sympos. Math. XXXIX}, 
Cambridge Univ. Press, Cambridge, 1999.
\nr
[BB4] M.T. Barlow, R.F. Bass. Divergence form operators on
fractal-like domains. {\it J. Funct. Analysis} {\bf 175} 
(2000), 214--247. 
\nr
[BB5] M.T. Barlow, R.F. Bass. Stability of parabolic Harnack inequalities. 
{\it Trans. Amer. Math. Soc.} {\bf 356} (2003) no. 4, 1501--1533.
\nr
[BBK] M.T. Barlow, R.F. Bass, T. Kumagai. 
Note on the equivalence of parabolic Harnack inequalities and heat kernel 
estimates. Unpublished note. 
\nl Available at: {\tt http://www.kurims.kyoto-u.ac.jp/\~{}kumagai/kumpre.html}
\nr
[BCG] M. Barlow, T. Coulhon, A. Grigor'yan. Manifolds and graphs with
slow heat kernel decay. {\it Invent. Math.} {\bf 144} (2001), 609-649. 
\nr 
[BCK] 
M. T. Barlow, T. Coulhon, T. Kumagai. 
Characterization of sub-Gaussian heat kernel 
estimates on strongly recurrent graphs. To appear {\it Comm. Pure Appl. Math.}
\nr 
[BP]
M.T. Barlow, E.A. Perkins. 
Brownian Motion on the Sierpi\'nski gasket. 
{\it Probab. Theory Rel. Fields \bf 79} (1988), 543--623.
\nr
[BPY] M. Barlow, J. Pitman, M. Yor. On Walsh's Brownian motions.  
{\it S\'em. Probab.} {\bf XXIII},  275--293, Lect. Notes Math. 
{\bf 1372}, Springer, Berlin, 1989. 
\nr
[BS] A. Bendikov, L. Saloff-Coste. 
On- and off-diagonal heat kernel behaviors on certain infinite dimensional 
Dirichlet spaces. {\it Amer. J. Math.} {\bf 122} (2000), 1205--1263. 
\nr
[BM] M. Biroli, U.A. Mosco. 
Saint-Venant type principle for Dirichlet forms on discontinuous media. 
{\it Ann. Mat. Pura Appl. \bf 169} (1995), 125--181. 
\nr
[Bom] E. Bombieri. Theory of minimal surfaces and a counter-example to the
Bernstein conjecture in high dimensions. Mimeographed Notes of Lectures held
at Courant Inst., New York Univ., 1970. 
\nr
[BG] E. Bombieri, E. Giusti. Harnack's inequality for elliptic
differential equations on minimal surfaces. 
{\it Invent. Math.} {\bf 15} (1972), 24--46. 
\nr
[Ch] Z.Q. Chen. On reflected Dirichlet spaces.  
{\it Probab. Theory Rel. Fields} {\bf 94} (1992), no. 2, 135--162.
\nr
[CS] T. Coulhon, L. Saloff-Coste. 
Vari\'et\'es riemanniennes isom\'etriques \`a l'infini. 
{\it  Rev. Math. Iberoamericana \bf 11} (1995), 687--726.
\nr
[Da] E.B. Davies. {\sl Heat kernels and spectral theory.} Cambridge
Univ. Press, 1989.
\nr
[Del] T. Delmotte. Parabolic Harnack inequality and estimates
of Markov chains on graphs.  {\it Rev. Math. Iberoamericana} {\bf 15}
(1999), 181--232. 
\nr
[FS] E.B. Fabes, D.W. Stroock. A new proof of Moser's parabolic
Harnack inequality using the old ideas of Nash.
{\it Arch. Mech. Rat. Anal. \bf 96} (1986), 327--338.
\nr 
[FHK]
P.J. Fitzsimmons, B.M. Hambly, T. Kumagai.
Transition density estimates for Brownian motion on affine nested fractals. 
{\it Comm. Math. Phys. \bf 165} (1994), 595--620.
\nr
[FOT] M. Fukushima, Y. Oshima, M. Takeda.
{\sl Dirichlet Forms and Symmetric Markov Processes}.
de Gruyter, Berlin, 1994.
\nr
[Gr1] A. Grigor'yan. The heat equation on noncompact Riemannian 
manifolds. 
(in Russian) {\it Matem. Sbornik. \bf 182} (1991), 55--87.
(English transl.) {\it Math. USSR Sbornik \bf 72} (1992), 47--77.
\nr 
[Gr2]  A. Grigor'yan. Heat kernels and function theory on metric measure
 spaces. In: {\sl Heat kernels and analysis on manifolds, graphs, and metric spaces} 
(Paris, 2002),  143--172, {\it Contemp. Math.} {\bf 338}, Amer. Math. Soc., 
Providence, RI, 2003.
\nr 
[GHL] 
A. Grigor'yan, J. Hu, K-S Lau.
Heat kernels on metric-measure spaces and an application to semi-linear 
elliptic equations. 
{\it Trans. Amer. Math. Soc.} {\bf 355} (2003) no.5, 2065--2095.
\nr
[GT1] A. Grigor'yan, A. Telcs. Sub-Gaussian estimates of heat kernels
on infinite graphs. {\it Duke Math. J.} {\bf 109} (2001), 452--510.
\nr 
[GT2] A. Grigor'yan, A. Telcs. Harnack inequalities and sub-Gaussian 
estimates for random walks.  {\it Math. Annalen} {\bf 324} (2002), 521--556. 
\nr 
[GT3] 
A. Grigor'yan, A. Telcs. 
Two-sided estimates of heat kernels in metric measure spaces. 
In preparation.
\nr
[Ha] P. Haj\l asz. Sobolev spaces on metric-measure spaces.  
In: {\it Heat kernels and analysis on manifolds, graphs, 
and metric spaces} (Paris, 2002), 173--218,  
{\it Contemp. Math. \bf 338}, Amer. Math. Soc., Providence, RI, 2003.
\nr 
[HaKo] P. Haj\l asz, P. Koskela. Sobolev met Poincar\'e. {\it Mem. Amer.
Math. Soc.} {\bf 145} (2000).
\nr 
[HK] B.M. Hambly, T. Kumagai. 
Transition density estimates for diffusion processes on post 
critically finite self-similar fractals.
{\sl Proc. London Math. Soc.} {\bf 78} (1999) no.3, 431--458.
\nr
[HK2]  B.M. Hambly, T. Kumagai. Heat kernel estimates for symmetric
random walks on a class of fractal graphs and stability under rough 
isometries. {\it Fractal geometry and applications: A Jubilee of 
B. Mandelbrot} (San Diego, CA, 2002), 233--260, {\sl Proc. Sympos. Pure Math.} 
{\bf 72}, Part 2, Amer. Math. Soc., Providence, RI, 2004.
\nr
[HKKZ] B.M. Hambly, T. Kumagai, S. Kusuoka, X.Y. Zhou. 
Transition density estimates for diffusion processes
on homogeneous random Sierpinski carpets. 
{\it J. Math. Soc. Japan} {\bf 52} (2000), 373--408.
\nr 
[HSC] W. Hebisch, L. Saloff-Coste. On the relation between
elliptic and parabolic Harnack inequalities. 
{\it Ann. Inst. Fourier (Grenoble)} {\bf  51}  (2001), 1437--1481. 
\nr 
[Hin]
M. Hino. On the singularity of energy measures on self-similar
sets. 
To appear {\it Probab. Theory Rel. Fields}. 
\nr 
[Je] 
D. Jerison. The weighted Poincar\'e inequality for vector
fields satisfying H\"ormander's condition.
{\it Duke Math. J.} {\bf 53} (1986), 503--523.
\nr 
[Jo] O.D. Jones. Transition probabilities for the simple random
walk on the Sierpinski graph. {\it Stoch. Proc. Appl. \bf 61} (1996),  45--69.
\nr
[Kan1] 
M. Kanai. Rough isometries and combinatorial approximations
of geometries of non-compact riemannian manifolds. {\it J. Math. Soc.
Japan} {\bf 37} (1985), 391--413. 
\nr 
[Kan2] 
M. Kanai. Analytic inequalities, and rough  isometries 
between  non-compact \break riemannian manifolds. 
In: {\it Curvature and topology of
Riemannian manifolds (Katata, 1985)}, 122--137, Lect. Notes Math. 
{\bf 1201}, Springer,  Berlin, 1986.
\nr
[Kas] A. Kasue.  
Convergence and energy forms in measure metric spaces.  
{\it S\=ugaku \bf 55} (2003), 20--36.
\nr
[Ke] S. Keith.  Modulus and the Poincar\'e inequality on metric measure 
spaces.  {\it Math. Z. \bf 245} (2003), 255--292.
\nr 
[Ki] J. Kigami. {\sl Analysis on fractals.}
Cambridge Tracts in Mathematics, {\bf 143},
Cambridge Univ. Press, Cambridge, 2001. 
\nr
[Kum1]
T. Kumagai.
Estimates of transition densities for Brownian motion on nested fractals. 
{\it Probab. Theory Rel. Fields \bf 96} (1993), 205--224.
\nr
[Kum2] T. Kumagai. Brownian motion penetrating fractals
-An application of the trace theorem of Besov spaces-. 
{\it J. Funct. Anal.} {\bf 170} (2000), 69--92.
\nr
[Kum3] 
T. Kumagai. Heat kernel estimates and parabolic Harnack inequalities
on graphs and resistance forms. {\sl Publ. RIMS Kyoto Univ. \bf 40} 
(2004), 793--818. 
\nr
[Kus] S. Kusuoka. Dirichlet forms on fractals and products of random
matrices. {\it Publ. RIMS Kyoto Univ.} {\bf 25} (1989), 659--680.
\nr 
[KS] S. Kusuoka, D.W. Stroock.  
Applications of the Malliavin calculus. III.  
{\it J. Fac. Sci. Univ. Tokyo Sect. IA Math. \bf 34} (1987), 391--442. 
\nr 
[KZ]
S. Kusuoka, X.Y. Zhou. Dirichlet form on fractals: Poincar\'e 
constant and resistance. {\it Probab. Theory Rel. Fields} 
{\bf 93} (1992), 169--196.
\nr 
[LJ]
Y. Le Jan. 
Mesures associees a une forme de Dirichlet. Applications. 
{\it Bull. Soc.  Math. France \bf 106} (1978), no. 1, 61--112.
\nr
[LiY] 
P. Li, S.-T. Yau. On the parabolic kernel of the Schr\"odinger 
operator. {\sl Acta Math.} {\bf 156} (1986) no. 3-4, 153--201.
\nr
[Mos]
U. Mosco. Composite media and asymptotic Dirichlet forms. 
{\it J. Funct. Anal. \bf 123} (1994), no. 2, 368--421.
\nr
[M1]  J. Moser. On Harnack's inequality for elliptic differential equations.
{\it Comm. Pure Appl. Math. \bf 14} (1961), 577--591.
\nr 
[M2] J. Moser. On Harnack's inequality for parabolic differential equations.
{\it Comm. Pure Appl. Math. \bf 17} (1964), 101--134. 
\nr
[M3] J. Moser. On a pointwise estimate for parabolic differential equations.
{\it Comm. Pure Appl. Math. \bf 24} (1971), 727--740.
\nr 
[N] J. Nash. Continuity of solutions of parabolic and elliptic equations. 
{\it Amer, J. Math. \bf 80}  (1958), 931--954.
\nr 
[PS] A. Phillips, D. Sullivan. Geometry of leaves.
{\it Topology} {\bf 20} (1981), 209--218. 
\nr
[SC1] L. Saloff-Coste. A note on Poincar\'e, Sobolev, and Harnack
inequalities. {\it Inter. Math. Res. Notices} {\bf 2} (1992), 27--38.
\nr
[SC2]  
L. Saloff-Coste. {\sl Aspects of Sobolev-Type inequalities.}
Lond. Math. Soc. Lect. Notes {\bf 289}, Cambridge Univ. Press, 2002.
\nr 
[St1] K.-T. Sturm.
Analysis on local Dirichlet spaces -II. Gaussian upper bounds for 
the fundamental solutions of parabolic Harnack equations.  
{\it Osaka J. Math. \bf 32} (1995), 275--312.
\nr
[St2] K.-T. Sturm. Analysis on local Dirichlet spaces -III. The parabolic
Harnack inequality. {\it J. Math. Pures Appl.} {\bf 75} (1996), no. 9,
273-297.
\nr 
[T] A. Telcs. Random walks on graphs with volume and time doubling.
Preprint.

\bs
\bs

\ni MTB: Department of Mathematics, University of British Columbia, Vancouver
V6T 1Z2, Canada
\ms
\ni RFB: Department of Mathematics, University of Connecticut, Storrs, 
CT 06269-3009, USA
\ms 
\ni TK: Research Institute for Mathematical 
Sciences, Kyoto University, Kyoto 606-8502, Japan

\bye